\documentclass{amsart}
\usepackage{amssymb}

\usepackage{amscd}

\newtheorem{theorem}{Theorem}[section]
\newtheorem{corollary}[theorem]{Corollary}

\newtheorem{lemma}[theorem]{Lemma}
\newtheorem{proposition}[theorem]{Proposition}

\theoremstyle{definition}
\newtheorem{definition}{Definition}[section]

\theoremstyle{remark}

\begin{document}
\title{Cohomology and Obstructions I: Geometry of formal Kuranishi theory}
\author{Herbert Clemens}
\address{Mathematics Department, University of Utah}
\email{clemens@math.utah.edu}
\date{August, 2002}
\maketitle

\begin{abstract}
The principle ``ambient cohomology of a Kaehler manifold annihilates
obstructions'' has been known and exploited since pioneering work of
Kodaira. This paper extends and unifies many known results in two contexts,
abstract deformations of compact Kaehler manifolds and deformations of
submanifolds within a given deformation of the ambient manifold.
\end{abstract}

\section{Introduction}

This paper is an attempt to simplify, clarify, and extend some results about
the interaction between deformation theory of K\"{a}hler manifolds and
submanifolds and cohomology of the ambient manifold. It was occasioned by a
review of the basics of Kuranishi theory in \cite{C1} and by the author's
desire to reconstruct the results of preprints \cite{R1}-\cite{R4} in the
setting of classical theory (see, for example, \S 3 of \cite{GM}). This
paper has two sequels (\cite{C1} and \cite{C2}), both dealing with
applications to the case of (deformations of) $K$-trivial threefolds.

The basic idea in this paper is always the same, namely:

\begin{quote}
Let $M_{0}$ be a compact K\"{a}hler manifold. Since the variations of Hodge
structure $H^{*}\left( M_{0}\right) $ over Artinian subschemes of $\Delta
_{0}$ always extend, therefore the subspaces of $H^{q}\left( \Omega
_{M_{0}}^{p}\right) $ associated to a particular geometric deformation
problem must pair to zero with the obstruction group for that problem.
\end{quote}

The ``therefore'' in the above assertion is not obvious and is based on the
equivalence of two pieces of data:

\begin{enumerate}
\item  The Gauss-Manin connection comparing the (trivial) topological
deformation to the (non-trivial) deformation of Hodge structure.

\item  The Kuranishi data associated to a topological trivialization of a
deformation of complex structure.
\end{enumerate}

Thus the main thrust of this paper is establishing the precise link between
the Kuranishi data and the Gauss-Manin connection.

The main use of the fact that obstructions pair to zero with certain
cohomology classes is to reduce the size of the obstruction space for
particular deformation problems. There are two main cases, many aspects of
which have already been treated by other authors (e.g. \cite{B}, \cite{BF}, 
\cite{FM}, \cite{Ka}, \cite{R1}, \cite{R2}, \cite{R3}, \cite{R4}, \cite{Ti}, 
\cite{To}). So the purpose here is to clarify and refine what these authors
and others have pointed out, namely that certain natural pairings between
ambient Hodge classes and obstructions measure nothing more than the
obstructions to deforming Hodge structures and therefore must vanish by \cite
{D}. This vanishing then gives useful limitations on the size of the
obstruction space in question.

\begin{quotation}
Case One: The obstructions $Obs$ to deforming a compact K\"{a}hler manifold $%
M_{0}$ annihilate the cohomology of $M_{0}$, that is, 
\begin{equation*}
Obs\otimes H^{p,q}\left( M_{0}\right)
\end{equation*}
lies in the nullspace of the natural pairing 
\begin{equation*}
H^{2}\left( T_{M_{0}}\right) \otimes H^{p,q}\left( M_{0}\right) \rightarrow
H^{p-1,q+2}\left( M_{0}\right) .
\end{equation*}

Case Two: Given a deformation $M/\Delta $ of a compact K\"{a}hler manifold $%
M_{0}$, and given a compact submanifold $Y_{0}$ such that the
sub-Hodge-structure 
\begin{equation*}
K_{0}^{r}=\sum K_{0}^{p,q}=\ker \left( H^{r}\left( M_{0}\right) \rightarrow
H^{r}\left( Y_{0}\right) \right)
\end{equation*}
deforms over $\Delta $ for some $r$, then obstructions $Obs$ to deforming $%
Y_{0}$ over $\Delta $ annihilate the primitive $r$-th cohomology of $M_{0}$,
that is, 
\begin{equation*}
Obs\otimes K_{0}^{p,q}
\end{equation*}
lies in the nullspace of the natural pairing 
\begin{equation*}
H^{1}\left( N_{Y_{0}\backslash M_{0}}\right) \otimes K_{0}^{p,q}\rightarrow
H^{p-1,q+1}\left( Y_{0}\right) .
\end{equation*}
\end{quotation}

\noindent For example, an immediate corollary of Case One is the fact that
all K\"{a}hler manifolds with trivial canonical bundle are unobstructed.

A variant of Case Two for non-compact $Y_{0}$ gives the analogous conclusion
in a relative setting, at least in the case of \textit{curvilinear}
deformations. For example, let $\Delta $ denote the complex unit disk and
let 
\begin{equation*}
p^{\prime }:Y_{0}\rightarrow Y_{0}^{\prime }
\end{equation*}
be a holomorphic family of $q$-dimensional compact submanifolds of $M_{0}$
with deformation $Y_{S}/Y_{S}^{\prime }/\Delta _{S}$ of $Y_{0}/Y_{0}^{\prime
}$ over a subscheme $S\subseteq \Delta $ such that, for some extendable
rational sub-Hodge-structure 
\begin{equation*}
K_{0}^{r}\subseteq H^{r}\left( M_{0}\right) ,
\end{equation*}
the natural map 
\begin{equation*}
K_{S}^{p+q+1,q-1}\overset{Gauss-Manin}{\longrightarrow }K_{S}^{p+q,q}%
\rightarrow R^{q}p_{*}^{\prime }\left( \Omega _{Y_{S}}^{p+q}\right)
\rightarrow \Omega _{Y_{S}^{\prime }}^{p}
\end{equation*}
is zero. Then the obstructions 
\begin{equation*}
Obs\subseteq R^{1}p_{*}^{\prime }\left( N_{Y_{0}\backslash Y_{0}^{\prime
}\times M_{0}}\right)
\end{equation*}
to extending the family $Y_{S}/Y_{S}^{\prime }$ again are such that 
\begin{equation*}
Obs\otimes K_{0}^{p+q+1,q-1}
\end{equation*}
lie in the nullspace of the pairing 
\begin{equation*}
R^{1}p_{*}^{\prime }\left( N_{Y_{0}\backslash Y_{0}^{\prime }\times
M_{0}}\right) \otimes K_{0}^{p+q+1,q-1}\rightarrow R^{q}p_{*}^{\prime
}\left( \Omega _{Y_{0}}^{p+q}\right) \rightarrow \Omega _{Y_{0}^{\prime
}}^{p}.
\end{equation*}
In this last result, the necessity of restricting one's attention to
curvilinear deformations is probably significant; that is, it only seems to
hold if we are working in the the reduced normal cone rather than in the
(larger) the normal sheaf (see, for example, \cite{BF}). We give several new
applications of this last result in the companion paper ``\textit{Cohomology
and Obstructions II: Curves on Calabi-Yau threefolds}.''

Finally, we combine Cases One and Two into the appropriate more general
theorem in the context of obstructions to deforming the pair $\left(
Y_{0},M_{0}\right) .$

As mentioned above, all these results derive from the very close
relationship between a Kuranishi datum associated to a deformation $M/\Delta 
$ and the Gauss-Manin connection associated to the deformation. We begin by
identifying (integrable) Kuranishi data 
\begin{equation*}
\xi \in A^{0,1}\left( T_{M_{0}}\right) \otimes \Bbb{C}\left[ \left[ t\right]
\right]
\end{equation*}
with $C^{\infty }$-trivializations of deformations $M/\Delta $ for which
each transverse fiber of the trivialization is a holomorphic copy of $\Delta 
$ (with holomorphic multi-parameter $t$). Letting 
\begin{equation*}
B^{*}\left( M/\Delta \right) =\frac{A^{*}\left( M/\Delta \right) }{
differential\ ideal\ \left\{ \overline{t},\partial t\right\} },
\end{equation*}
we have a (formal) isomorphism of $d$-differential graded algebras 
\begin{equation*}
\varphi :B^{*}\left( M/\Delta \right) \rightarrow A^{*}\left( M_{0}\right)
\otimes \Bbb{C}\left[ \left[ t\right] \right]
\end{equation*}
for which 
\begin{equation*}
\varphi B^{p,q}\left( M/\Delta \right) =e^{\left\langle \left. \xi \right| \
\right\rangle }A^{p,q}\left( M_{0}\right) \otimes \Bbb{C}\left[ \left[
t\right] \right] .
\end{equation*}
(See the notational comment at the end of this introduction.)

Defining Lie differentiation 
\begin{eqnarray*}
L_{\xi } &=&\left[ \left\langle \left. \xi \right| \ \right\rangle ,d\right]
\\
L_{\xi }^{1,0} &=&\left[ \left\langle \left. \xi \right| \ \right\rangle
,\partial \right] ,
\end{eqnarray*}
we have the following table of correspondences 
\begin{equation}
\begin{array}{ccccc}
B^{p,q}\left( M/\Delta \right) & \leftrightarrow & \varphi B^{p,q}\left(
M/\Delta \right) & \overset{e^{\left\langle \left. \xi \right| \
\right\rangle }}{\leftrightarrows } & A^{p,q}\left( M_{0}\right) \otimes 
\Bbb{C}\left[ \left[ t\right] \right] \\ 
d & \leftrightarrow & d &  &  \\ 
\overline{\partial }_{M/\Delta } & \leftrightarrow & \overline{\partial }
_{M_{0}}-L_{\xi } & \leftrightarrow & \overline{\partial }_{M_{0}}-L_{\xi
}^{1,0}
\end{array}
\label{T1}
\end{equation}
(See \S \ref{newGM}.) These identifications are the key to the results in
this paper.

This paper has three parts. Part \ref{part1} attempts to bring together a
self-contained account of formal Kuranishi theory from the point of view
which will be used in the applications. It is included because the point of
view, namely ``calculus modulo $\overline{t}$,'' seems to be a bit different
from the usual one, and, I believe, makes the essentials of the theory more
transparent. At the cost of a few explicit local computations at the outset,
all standard results derive directly from the table $\left( \ref{T1}\right) $
and this framework easily provides several new applications. The
significance of one of the explicit local computations, the computation of
local gauge transformations in \S \ref{lg}, may not be apparent from the
applications of Kuranishi theory derived in this paper, so that section \S 
\ref{lg} may well be skipped for the purposes of those applications. However
I beg the reader's indulgence on this point. The section is included for
completeness, and its inclusion is perhaps additionally justified by its
central role in applications which appear elsewhere, such as in \cite{C3}.

Part \ref{part4} contains the main new results of this paper and their
proofs. Finally Part \ref{part5} is an Appendix containing a full proof of
the ``intuitively obvious'' fact that transversely holomorphic
trivializations always exist, together with a proof of the classical
characterization of integrability of Kuranishi data and of the standard Lie
derivative identities which are used in that proof and elsewhere in this
paper.

A word about notation: The holomorphic tangent space for a complex manifold $%
M_{0}$ will be denoted by $T_{M_{0}}$ or alternatively by $T_{1,0}\left(
M_{0}\right) $ depending on the context, with the full complexified tangent
space and cotangent spaces denoted respectively as 
\begin{eqnarray*}
T\left( M_{0}\right) &=&T_{1,0}\left( M_{0}\right) \oplus T_{0,1}\left(
M_{0}\right) \\
T^{*}\left( M_{0}\right) &=&T^{1,0}\left( M_{0}\right) \oplus T^{0,1}\left(
M_{0}\right) .
\end{eqnarray*}
We will write 
\begin{equation*}
T^{p,q}\left( M_{0}\right) =\bigwedge\nolimits^{p}T^{1,0}\left( M_{0}\right)
\wedge \bigwedge\nolimits^{q}T^{0,1}\left( M_{0}\right)
\end{equation*}
and will denote the space of $C^{\infty }$-forms of type $\left( p,q\right) $
with coefficients in, for example, $T_{M_{0}}$ as $A^{p,q}\left(
T_{M_{0}}\right) .$ Also, for an operator $D$, the notations 
\begin{equation*}
e^{D}
\end{equation*}
and 
\begin{equation*}
\exp D
\end{equation*}
will be used interchangeably, denoting the formal power series of operators 
\begin{equation*}
\sum\nolimits_{i\geq 0}\frac{1}{i!}\underset{i\ times}{\underbrace{D\circ
\ldots \circ D}}.
\end{equation*}

\part{Formal Kuranishi theory\label{part1}}

\section{Transversely holomorphic trivializations, local version\label{1}}

Let $D$ be a complex polydisk with holomorphic coordinate system $x=\left(
x_{i}\right) $ and let $\Delta $ be a complex polydisk with holomorphic
coordinate system $t=\left( t_{j}\right) $ centered at $t=0$. We use the
standard notation 
\begin{equation*}
t^{J}:=t_{1}^{j_{1}}\cdot \ldots \cdot t_{m}^{j_{m}}
\end{equation*}
whenever $J=\left( j_{1},\ldots ,j_{m}\right) $ is an $m$-tuple of
non-negative integers. Also 
\begin{equation*}
\left| J\right| :=\sum_{k=1}^{m}j_{k}.
\end{equation*}
Suppose I have a holomorphic map 
\begin{eqnarray*}
\Delta &\rightarrow &D \\
t &\mapsto &x+\varphi \left( t\right)
\end{eqnarray*}
with $\varphi \left( 0\right) =0$. Then, for a function $f$ on $D$, at least
formally the formula for the pull-back function is 
\begin{equation*}
f\circ \varphi =e^{L_{\beta +\overline{\beta }}}\left( f\right)
\end{equation*}
where 
\begin{equation*}
\beta =\sum\nolimits_{i}\varphi _{i}\left( t\right) \cdot \frac{\partial }{
\partial x_{i}}=\sum\nolimits_{i,J,\left| J\right| >0}a_{J,i}\cdot
t^{J}\cdot \frac{\partial }{\partial x_{i}}
\end{equation*}
is the vector field giving the flow $\varphi \left( t,s\right) =x+s\varphi
\left( t\right) $ and $L$ denotes Lie differentiation.

Next suppose that we have a $C^{\infty }$-diffeomorphism 
\begin{eqnarray*}
G &:&D\times \Delta \rightarrow D\times \Delta \\
\left( x,t\right) &\mapsto &\left( x+\varphi \left( x,t\right) ,t\right)
\end{eqnarray*}
such that $\varphi \left( x,t\right) $ is holomorphic in $t$ and 
\begin{equation*}
\varphi \left( x,0\right) =0.
\end{equation*}
Then, for fixed $x$ and a function $f$ on $D$, we consider $f$ as a function
on $D\times \Delta $ via the product structure and we have, again at least
formally, that 
\begin{equation*}
G^{*}\left( f\right) =e^{L_{\beta \left( x\right) +\overline{\beta \left(
x\right) }}}\left( f\right)
\end{equation*}
where 
\begin{equation*}
\beta \left( x\right) =\sum\nolimits_{J,\left| J\right| >0}a_{J,i}\left(
x\right) \frac{\partial }{\partial x_{i}}t^{J}
\end{equation*}
for some $C^{\infty }$-functions $a_{J,i}\left( x\right) $ on $D$. So now if 
\begin{equation*}
f\left( x,t\right) =\sum f_{J,K}\left( x\right) t^{J}\overline{t}^{K}
\end{equation*}
then the formula 
\begin{equation}
G^{*}\left( f\left( x,t\right) \right) =e^{L_{\beta +\overline{\beta }
}}\left( f\left( x,t\right) \right)  \label{n1.0}
\end{equation}
continues to hold (and this formula is true geometrically for real analytic
functions). We call the ring 
\begin{equation*}
\frak{C}=\left\{ \sum f_{J,K}t^{J}\overline{t}^{K}:f_{J,K}\in C^{\infty
}\left( D\right) \right\}
\end{equation*}
``the ring of $\Delta $-formal functions'' and form the quotient ring 
\begin{equation}
\frac{\frak{C}}{\left\{ \overline{t}\right\} }\cong \left( C^{\infty }\left(
D^{n}\right) \otimes \Bbb{C}\left[ \left[ t\right] \right] \right)
\label{n1.1}
\end{equation}
by dividing by the ideal generated by the functions $\overline{t}_{j}$. Thus:

\begin{proposition}
\label{p1.1}Formally 
\begin{equation*}
G^{*}=e^{L_{\beta }}:\frac{\frak{C}}{\left\{ \overline{t}\right\} }%
\rightarrow \frac{\frak{C}}{\left\{ \overline{t}\right\} }.
\end{equation*}
\end{proposition}

\begin{proof}
Immediate from $\left( \ref{n1.0}\right) $ and $\left( \ref{n1.1}\right) $.
\end{proof}

We want to use the diffeomorphism 
\begin{equation}
G:D\times \Delta \rightarrow D\times \Delta  \label{1.5}
\end{equation}
to induce a new complex structure on the domain. If we denote 
\begin{equation}
F_{\beta }=G^{-1},  \label{n1.3}
\end{equation}
then 
\begin{equation*}
F_{\beta }^{*}=e^{-L_{\beta }}:\frac{\frak{C}}{\left\{ \overline{t}\right\} }
\rightarrow \frac{\frak{C}}{\left\{ \overline{t}\right\} }
\end{equation*}
and it makes sense to ask for the 
\begin{equation*}
f\in \left( C^{\infty }\left( D\right) \otimes \Bbb{C}\left[ \left[ t\right]
\right] \right)
\end{equation*}
such that 
\begin{equation*}
F_{\beta }^{*}\left( f\right)
\end{equation*}
is holomorphic, that is, 
\begin{equation*}
\left( \overline{\partial }\circ e^{-L_{\beta }}\right) \left( f\right) =0
\end{equation*}
or, what is the same, 
\begin{equation*}
\left( e^{L_{\beta }}\circ \overline{\partial }\circ e^{-L_{\beta }}\right)
\left( f\right) =0.
\end{equation*}
Referring to the definition in $\left( \ref{2.7.1}\right) $ below, we
rewrite this as 
\begin{equation}
\left( \overline{\partial }-L_{\xi _{\beta }}\right) \left( f\right) =0
\label{1.7}
\end{equation}
where, as operators on $\frac{\frak{C}}{\left\{ \overline{t}\right\} }$, 
\begin{equation}
L_{\xi _{\beta }}=\left[ \overline{\partial },e^{L_{\beta }}\right] \circ
e^{-L_{\beta }}.  \label{1.4}
\end{equation}
Since $L_{\xi _{\beta }}$ is a $\left( 0,1\right) $-form with coeffiecients
which are first-order which annihilate constant functions, the power series 
\begin{equation}
\xi _{\beta }=\sum\nolimits_{\left| J\right| >0}\xi _{\beta ,J}t^{J}
\label{1.3}
\end{equation}
has coefficients $\xi _{\beta ,J}$ which are $\left( 0,1\right) $-forms in $%
T_{1,0}\left( D\right) $. $\xi _{\beta }$ is called the \textit{Kuranishi
data} for the local gauge transformation $F_{\beta }$. Notice that 
\begin{equation*}
\left( \overline{\partial }-L_{\xi _{\beta }}\right) :\frac{\frak{C}}{%
\left\{ \overline{t}\right\} }\rightarrow \frac{\frak{C}}{\left\{ \overline{t%
}\right\} }\otimes _{C^{\infty }\left( D\right) }A^{0,1}\left( D\right)
\end{equation*}
or, what is the same 
\begin{equation*}
\left( \overline{\partial }-L_{\xi _{\beta }}\right) :A^{0}\left( D\right)
\otimes \Bbb{C}\left[ \left[ t\right] \right] \rightarrow A^{0,1}\left(
D\right) \otimes \Bbb{C}\left[ \left[ t\right] \right] .
\end{equation*}

\section{One-parameter gauge transformations, local version\label{lg}}

Consider now a family 
\begin{equation*}
\beta \left( s\right) \in A^{0}\left( T_{1,0}\left( D\right) \right) \otimes 
\Bbb{C}\left[ \left[ t\right] \right]
\end{equation*}
which is holomorphic in $s$ and $t$. Correspondingly we have a one complex
parameter family 
\begin{eqnarray*}
F_{\beta \left( s\right) }^{-1} &:&D\times \Delta \rightarrow D\times \Delta
\\
\left( x,t\right) &\mapsto &\left( x+\varphi \left( x,t,s\right) ,t\right)
\end{eqnarray*}
with the property that, for fixed $x\in D$, $F_{\beta \left( s\right) }^{-1}$
is holomorphic in $s$ and $t$. Set 
\begin{equation*}
\gamma \left( s\right) =\frac{\partial \beta \left( s\right) }{\partial s}
=\sum\nolimits_{i}\frac{\partial \varphi _{i}}{\partial s}\frac{\partial }{
\partial x_{i}}\in T_{1,0}\left( D\right) \otimes \Bbb{C}\left[ \left[
t\right] \right]
\end{equation*}
which is therefore analytic in $s$. Then 
\begin{eqnarray*}
\left. \frac{\partial f\left( x+\varphi \left( x,t,s^{\prime }\right)
,t\right) }{\partial s^{\prime }}\right| _{s\prime =s} &=&\left( L_{\gamma
\left( s\right) }f\right) \left( x+\varphi \left( x,t,s\right) ,t\right) \\
&=&\left( F_{\beta \left( s\right) }^{-1}\right) ^{*}\left( L_{\gamma \left(
s\right) }f\right)
\end{eqnarray*}
so 
\begin{equation*}
\frac{\partial \left( F_{\beta \left( s\right) }^{-1}\right) ^{*}\left(
f\right) }{\partial s}=\left( F_{\beta \left( s\right) }^{-1}\right)
^{*}\left( L_{\gamma \left( s\right) }f\right) .
\end{equation*}
So as operators on 
\begin{equation*}
\frac{\frak{C}}{\left\{ \overline{t}\right\} }
\end{equation*}
we have 
\begin{equation*}
\frac{\partial e^{L_{\beta \left( s\right) }}}{\partial s}=e^{L_{\beta
\left( s\right) }}\circ L_{\gamma \left( s\right) }.
\end{equation*}

Now 
\begin{eqnarray}
\left( F_{\beta \left( s\right) }\right) _{*}\left( L_{\gamma \left(
s\right) }\right) &=&\left( F_{\beta \left( s\right) }^{-1}\right) ^{*}\circ
L_{\gamma \left( s\right) }\circ F_{\beta \left( s\right) }^{*}
\label{pushforward} \\
&=&L_{\tilde{\alpha}\left( s\right) }  \notag
\end{eqnarray}
for 
\begin{equation*}
\tilde{\alpha}\left( s\right) =\left( F_{\beta \left( s\right) }\right)
_{*}\left( \gamma \left( s\right) \right) \in \left( T_{1,0}\left( D\right)
\oplus T_{0,1}\left( D\right) \right) \otimes \Bbb{C}\left[ \left[ t,%
\overline{t}\right] \right] .
\end{equation*}
And so as operators on $\frac{\frak{C}}{\left\{ \overline{t}\right\} }$, $%
\left( \ref{pushforward}\right) $ becomes the equality 
\begin{equation*}
e^{L_{\beta \left( s\right) }}\circ L_{\gamma \left( s\right) }\circ
e^{-L_{\beta \left( s\right) }}=L_{\alpha \left( s\right) }
\end{equation*}
and the vector field $\alpha \left( s\right) $ is of type $\left( 1,0\right) 
$ since $\beta \left( s\right) $ is. As an immediate corollary of the
definition of $\alpha $ we have 
\begin{equation*}
\left. \frac{\partial \left( F_{\beta \left( s^{\prime }\right)
}^{-1}\right) ^{*}}{\partial s^{\prime }}\right| _{s^{\prime }=s}=L_{\alpha
\left( s\right) }:\frac{\frak{C}}{\left\{ \overline{t}\right\} }\rightarrow 
\frac{\frak{C}}{\left\{ \overline{t}\right\} }.
\end{equation*}

\begin{lemma}
\label{l1.1} 
\begin{equation*}
\frac{\partial \xi _{\beta \left( s\right) }}{\partial s}=\overline{\partial 
}\alpha \left( s\right) +\left[ \alpha \left( s\right) ,\xi _{\beta \left(
s\right) }\right] .
\end{equation*}
\end{lemma}

\begin{proof}
First we compute 
\begin{eqnarray*}
\frac{\partial \left( e^{-L_{\beta \left( s\right) }}\right) }{\partial s}
&=&-e^{-L_{\beta \left( s\right) }}\circ \frac{\partial e^{L_{\beta \left(
s\right) }}}{\partial s}\circ e^{-L_{\beta \left( s\right) }} \\
&=&-L_{\gamma \left( s\right) }\circ e^{-L_{\beta \left( s\right) }}.
\end{eqnarray*}

Then 
\begin{eqnarray*}
\frac{\partial \left( \overline{\partial }-L_{\xi _{\beta \left( s\right)
}}\right) }{\partial s} &=&\frac{\partial \left( e^{L_{\beta \left( s\right)
}}\circ \overline{\partial }\circ e^{-L_{\beta \left( s\right) }}\right) }{
\partial s} \\
&=&\frac{\partial \left( e^{L_{\beta \left( s\right) }}\right) }{\partial s}
\circ \overline{\partial }\circ e^{-L_{\beta \left( s\right) }}+e^{L_{\beta
\left( s\right) }}\circ \overline{\partial }\circ \frac{\partial \left(
e^{-L_{\beta \left( s\right) }}\right) }{\partial s} \\
&=&e^{L_{\beta \left( s\right) }}\circ L_{\gamma \left( s\right) }\circ 
\overline{\partial }\circ e^{-L_{\beta \left( s\right) }}-e^{L_{\beta \left(
s\right) }}\circ \overline{\partial }\circ L_{\gamma \left( s\right) }\circ
e^{-L_{\beta \left( s\right) }} \\
&=&e^{L_{\beta \left( s\right) }}\circ \left[ L_{\gamma \left( s\right) },%
\overline{\partial }\right] \circ e^{-L_{\beta \left( s\right) }} \\
&=&\left[ e^{L_{\beta \left( s\right) }}\circ L_{\gamma \left( s\right)
}\circ e^{-L_{\beta \left( s\right) }},\overline{\partial }-L_{\xi _{\beta
\left( s\right) }}\right] \\
&=&\left[ L_{\alpha \left( s\right) },\overline{\partial }-L_{\xi _{\beta
\left( s\right) }}\right] \\
&=&-L_{\overline{\partial }\alpha \left( s\right) }-L_{\left[ \alpha \left(
s\right) ,\xi _{\beta \left( s\right) }\right] }
\end{eqnarray*}
where the last equality anticipates the definition in $\left( \ref{2.7.2}
\right) $.
\end{proof}

In the other direction, given two $t$-holomorphic maps 
\begin{eqnarray*}
D\times \Delta &\rightarrow &D\times \Delta \\
t &\mapsto &\left( x+\varphi ^{0}\left( x,t\right) ,t\right) \\
t &\mapsto &\left( x+\varphi ^{1}\left( x,t\right) ,t\right)
\end{eqnarray*}
with $\varphi ^{0}\left( x,0\right) =\varphi ^{1}\left( x,0\right) =0$, we
consider the family of maps 
\begin{equation*}
\left( x+\varphi ^{0}\left( x,t\right) +s\left( \varphi ^{1}\left(
x,t\right) -\varphi ^{0}\left( x,t\right) \right) ,t\right)
\end{equation*}
joining them. Then 
\begin{equation*}
\gamma \left( s\right) =\gamma =\sum\nolimits_{i}\left( \varphi
_{i}^{1}-\varphi _{i}^{0}\right) \frac{\partial }{\partial x_{i}}
\end{equation*}
is independent of $s$. As operators on $\frac{\frak{C}}{\left\{ \overline{t}
\right\} }$, we write 
\begin{equation*}
e^{L_{\beta \left( s\right) }}=e^{L_{\beta \left( 0\right) }}\circ
e^{L_{s\gamma }}
\end{equation*}
so that 
\begin{eqnarray*}
e^{L_{\beta \left( s\right) }}\circ L_{\gamma \left( s\right) }\circ
e^{-L_{\beta \left( s\right) }} &=&e^{L_{\beta \left( 0\right) }}\circ
e^{L_{s\gamma }}\circ L_{\gamma }\circ e^{-L_{s\gamma }}\circ e^{-L_{\beta
\left( 0\right) }} \\
&=&e^{L_{\beta \left( 0\right) }}\circ L_{\gamma }\circ e^{-L_{\beta \left(
0\right) }} \\
&=&L_{\alpha }
\end{eqnarray*}
where $\alpha $ is also independent of $s$.

\begin{lemma}
\label{l1.2} For the analytic family $\beta \left( s\right) $ for which $%
\beta \left( 0\right) =\beta $ and 
\begin{equation*}
\beta \left( s\right) =\beta +s\gamma 
\end{equation*}
we have 
\begin{equation*}
\xi _{\beta \left( s\right) }=\exp \left( \left[ s\alpha ,\ \right] \right)
\left( \xi _{\beta }\right) +\frac{\exp \left( \left[ s\alpha ,\ \right]
\right) -1}{\left[ s\alpha ,\ \right] }\left( s\left[ \bar{\partial },\alpha
\right] \right) .
\end{equation*}
\end{lemma}

\begin{proof}
Call the right-hand-side $\xi _{s}$. We compute the coefficient of $h$ in
the power series expansion of 
\begin{equation*}
\xi _{s+h}=\exp \left( \left[ \left( s+h\right) \alpha ,\ \right] \right)
\left( \xi _{\beta }\right) -\frac{1-\exp \left( \left[ \left( s+h\right)
\alpha ,\ \right] \right) }{\left[ \left( s+h\right) \alpha ,\ \right] }
\left( \left( s+h\right) \left[ \bar{\partial },\alpha \right] \right) ,
\end{equation*}
which, by a straightforward computation is 
\begin{equation*}
\left[ \alpha ,\ \right] \exp \left( \left[ s\alpha ,\ \right] \right)
\left( \xi _{\beta }\right) +\sum\nolimits_{k=0}^{\infty }\frac{\left[
s\alpha ,\ \right] ^{k}}{k!}\left[ \bar{\partial },\alpha \right]
\end{equation*}
or 
\begin{equation*}
\left[ \alpha ,\ \right] \left( \exp \left( \left[ s\alpha ,\ \right]
\right) \left( \xi _{\beta }\right) +\sum\nolimits_{k=1}^{\infty }\frac{%
\left[ s\alpha ,\ \right] ^{k-1}}{k!}\left[ \bar{\partial },\alpha \right]
\right) +\left[ \bar{\partial },\alpha \right]
\end{equation*}
and so is computed as 
\begin{equation*}
\left[ \alpha ,\xi _{s}\right] +\left[ \bar{\partial },\alpha \right] .
\end{equation*}
So by Lemma \ref{l1.1} 
\begin{equation*}
\frac{\partial \left( \xi _{\beta \left( s\right) }-\xi _{s}\right) }{
\partial s}=\left[ \alpha ,\left( \xi _{\beta \left( s\right) }-\xi
_{s}\right) \right]
\end{equation*}
and 
\begin{equation*}
\xi _{\beta \left( 0\right) }-\xi _{0}=0.
\end{equation*}
Since $\xi _{\beta \left( s\right) }-\xi _{s}$ is analytic in $s$, we
conclude 
\begin{equation*}
\xi _{\beta \left( s\right) }=\xi _{s}.
\end{equation*}
\end{proof}

\section{Transversely holomorphic trivializations\label{2}}

We next interpret the Newlander-Nirenberg-Kuranishi theory of deformations
of complex structures (see \cite{Ku}, Chapter 5 of \cite{Ko}, II.1 of \cite
{G} and \cite{GM}) in terms of ``transversely holomorphic'' trivializations
of a deformation. Let 
\begin{equation}
M\overset{\pi }{\longrightarrow }\Delta =\left\{ t=\left( t_{1},\ldots
,t_{m}\right) :\left| t_{j}\right| <<1,\forall j\right\} .  \label{2.1}
\end{equation}
be a deformation of a compact complex manifold K\"{a}hler manifold $M_{0}$
of dimension $m$.

\begin{definition}
\label{2.2}A $C^{\infty }$-projection 
\begin{equation*}
M\overset{\sigma }{\longrightarrow }M_{0}
\end{equation*}
will be called transversely holomorphic if all its fibers are complex
holomorphic disks meeting $M_{0}$ transversely. If $\sigma $ is a
transversely holomorphic projection, the diffeomorphism 
\begin{equation*}
F_{\sigma }:M\overset{\left( \sigma ,\pi \right) }{\longrightarrow }M_{0}%
\times \Delta 
\end{equation*}
will be called transversely holomorphic trivialization.
\end{definition}

\begin{proposition}
\label{p2.3}A transversely holomorphic trivialization is a diffeomorphism 
\begin{equation*}
F_{\sigma }:M\overset{\left( \sigma ,\pi \right) }{\longrightarrow }M_{0}%
\times \Delta 
\end{equation*}
for which there is a covering of $M$ by analytic open sets $W$ such that 
\begin{equation*}
\left. F_{\sigma }\right| _{W}=\left. F_{\beta }\right| _{W}
\end{equation*}
for some diffeomorphism $G$ as in $\left( \ref{1.5}\right) $ and $F_{\beta }$
as in $\left( \ref{n1.3}\right) $.
\end{proposition}

It is intuitively ``clear'' that transversely holomorphic trivializations
exist for any deformation $M/\Delta $. After all, these are just $C^{\infty
} $-trivializations for which the deformation trajectory of each fixed point 
$x_{0}\in M_{0}$ is holomorphic. We are just fitting those holomorphic
trajectories together in a $C^{\infty }$ way as $x_{0}$ moves on $M_{0}.$ A
precise proof of the existence of transversely holomorphic trivializations
and their properties is given in \S \ref{A1} of the Appendix to this paper.

Given any transversely holomorphic trivialization $F_{\sigma }$, under the $%
C^{\infty }$-isomorphisms 
\begin{equation*}
M_{t}\cong M_{0}
\end{equation*}
induced by $\sigma $, the holomorphic cotangent space of $M_{t}$ corresponds
to a subspace 
\begin{equation*}
T_{t}^{1,0}\subseteq T_{M_{0}}^{*}.
\end{equation*}
If 
\begin{equation*}
\pi ^{1,0}+\pi ^{0,1}:T_{M_{0}}^{*}\rightarrow T_{M_{0}}^{1,0}\oplus
T_{M_{0}}^{0,1}
\end{equation*}
are the two projections, the projection 
\begin{equation*}
\pi ^{1,0}:T_{t}^{1,0}\rightarrow T_{M_{0}}^{1,0}
\end{equation*}
is an isomorphism for small $t$ so that the composition 
\begin{equation*}
T_{M_{0}}^{1,0}\overset{\left( \pi ^{1,0}\right) ^{-1}}{\longrightarrow }
T_{t}^{1,0}\overset{\pi ^{0,1}}{\longrightarrow }T_{M_{0}}^{0,1},
\end{equation*}
gives a $C^{\infty }$-mapping 
\begin{equation*}
\xi \left( t\right) :T_{M_{0}}^{1,0}\rightarrow T_{M_{0}}^{0,1}
\end{equation*}
which determines and is determined by the deformation of (almost) complex
structure. Said differently, $T_{t}^{1,0}$ is uniquely determined as the
graph of $\xi \left( t\right) $.

Furthermore, by \S \ref{1} and Proposition \ref{p2.3}, for each commutative
diagram 
\begin{equation*}
\begin{array}{ccc}
D\times \Delta & \overset{i}{\hookrightarrow } & M \\ 
\downarrow &  &  \\ 
\Delta & = & \Delta
\end{array}
,
\end{equation*}
we have 
\begin{equation*}
i^{*}\xi =\xi _{\beta }
\end{equation*}
for $\xi _{\beta }$ as in $\left( \ref{1.4}\right) $. By this fact, or by
the explicit local formula given in Lemma \ref{2.3} of \S \ref{A1}$,$%
\begin{equation*}
\xi \left( t\right) =\sum\nolimits_{J,\left| J\right| >0}\xi _{J}t^{J}
\end{equation*}
is holomorphic in $t$.

\begin{lemma}
\label{2.5}The transversely holomorphic trivialization $\sigma $
distinguishes a subset of the complex-valued $C^{\infty }$-functions on $M$,
namely those $C^{\infty }$-functions $f$ which restrict to a holomorphic
function on each fiber of $\sigma $ (i.e. functions with a power-series
representation 
\begin{equation}
\sum\nolimits_{I}f_{I}t^{I}:=\sum\nolimits_{I}\left( f_{I}\circ \sigma
\right) t^{I}  \label{2.5.1}
\end{equation}
where the $f_{I}$ are $C^{\infty }$-functions on $M_{0}$). Furthermore, for $%
\left| J\right| >0$, there exist elements 
\begin{equation*}
\xi _{J}\in A^{0,1}\left( T_{M_{0}}^{1,0}\right) 
\end{equation*}
such that a function $\left( \ref{2.5.1}\right) $ is holomorphic if and only
if: 
\begin{equation}
\begin{tabular}{l}
$\bar{D}_{\sigma }\left( f\right) :=\left( \bar{\partial }_{M_{0}}-\xi
\right) \left( f\right) =$ \\ 
$\sum\nolimits_{I}\bar{\partial }f_{I}t^{I}-\sum\nolimits_{I,J,\left|
J\right| >0}\xi _{J}\left( f_{I}\right) t^{I+J}=0.$%
\end{tabular}
\label{2.5.2}
\end{equation}
Let $\sigma _{t}$ denote the restriction of $\sigma $ to $M_{t}$. If $\left(
v^{l}\right) $ are local holomorphic coordinates on $M_{0}$, a local basis
for 
\begin{equation*}
\left( \sigma _{t}^{-1}\right) ^{*}\left( T_{M_{t}}^{1,0}\right) 
\end{equation*}
is given by 
\begin{equation}
\alpha ^{l}:=dv^{l}+\left\langle \left. \xi \right| dv^{l}\right\rangle .
\label{2.5.3}
\end{equation}
Thus, modulo $\overline{t}$, the subspace 
\begin{equation*}
\left( \sigma _{t}^{-1}\right) ^{*}\left( T_{M_{t}}^{0,1}\right) 
\end{equation*}
is ``stationary,'' that is 
\begin{equation*}
\left( \sigma _{t}^{-1}\right) ^{*}\left( T_{M_{t}}^{0,1}\right) \equiv
T_{M_{0}}^{0,1}.
\end{equation*}
\end{lemma}

\begin{proof}
The formulas $\left( \ref{2.5.1}\right) $-($\left( \ref{2.5.2}\right) $
follow directly from $\left( \ref{1.4}\right) $ or from the formulas in
Lemma \ref{2.3} of \S \ref{A1} in the Appendix. For $\left( \ref{2.5.3}
\right) $, see Lemma \ref{2.3}.
\end{proof}

Also by $\left( \ref{1.4}\right) $ or the uniqueness assertion \ref{2.3.3}
of Lemma \ref{2.3} we have:

\begin{proposition}
\label{2.6}Two deformation/trivializations as in $\left( \ref{2.1}\right) $
and Definition \ref{2.2} are related by a holomorphic isomorphism over $%
\Delta ,$ that is, there is a commutative diagram 
\begin{equation*}
\begin{tabular}{ccc}
$M$ & $\overset{\varphi }{\longrightarrow }$ & $M^{\prime }$ \\ 
$\downarrow \sigma $ &  & $\downarrow \sigma ^{\prime }$ \\ 
$M_{0}$ & $=$ & $M_{0},$%
\end{tabular}
\end{equation*}
if and only if 
\begin{equation*}
\bar{D}_{\sigma }=\bar{D}_{\sigma ^{\prime }}.
\end{equation*}
\end{proposition}

\begin{proof}
One implication is immediate from the definitions of $\bar{D}_{\sigma }$ and 
$\bar{D}_{\sigma ^{\prime }}$. For the other, by Lemma \ref{2.3} in \S \ref
{A1} the equality 
\begin{equation*}
\xi _{\sigma }=\xi _{\sigma ^{\prime }}
\end{equation*}
implies that the codifferential of the $C^{\infty }$-automorphism 
\begin{equation*}
\varphi :=\left( \sigma ^{\prime },\pi \right) ^{-1}\circ \left( \sigma ,\pi
\right) :M\rightarrow M^{\prime }
\end{equation*}
preserves the $\left( 1,0\right) $-subspace of the (complexified) cotangent
space and therefore $\varphi $ is holomorphic. That is, $T_{t}^{1,0}$ above
is uniquely determined as the graph of $\xi \left( t\right) $ so that 
\begin{eqnarray*}
\left( \left( \sigma ^{\prime },\pi \right) ^{-1}\right) ^{*}\left(
T_{M^{\prime }/\Delta }^{1,0}\right) &=&T_{t}^{1,0} \\
\left( \sigma ,\pi \right) ^{*}\left( \left( \left( \sigma ^{\prime },\pi
\right) ^{-1}\right) ^{*}\left( T_{M^{\prime }/\Delta }^{1,0}\right) \right)
&=&T_{M/\Delta }^{1,0}
\end{eqnarray*}
and so $\varphi $ is holomorphic.
\end{proof}

If, for an ideal $\frak{A}\subseteq \frak{m}=\left\{ t_{1},\ldots
,t_{s}\right\} $, we let 
\begin{equation*}
\Delta _{\frak{A}}:=Spec\frac{\Bbb{C}\left[ t\right] }{\frak{A}}\subseteq
\Delta .
\end{equation*}
Then using the partition-of-unity $\left\{ \rho _{W}\right\} $ one easily
sees that any holomorphic function on 
\begin{equation*}
M_{\frak{A}}:=\pi ^{-1}\left( \Delta _{\frak{A}}\right)
\end{equation*}
extends to a function in the class $\left( \ref{2.5.1}\right) $. Also, if 
\begin{equation*}
f=\sum\nolimits_{I}f_{I}t^{I}
\end{equation*}
restricts to a holomorphic function on $M_{\frak{A}}$, for $\left| I\right|
\leq n$ we have 
\begin{equation*}
\overline{\partial }-\xi \left( f\right) \subseteq A^{0,1}\left(
M_{0}\right) \otimes \frak{A\subseteq }A^{0,1}\left( M_{0}\right) \otimes 
\Bbb{C}\left[ t\right] .
\end{equation*}

\section{Gauge transformations}

An immediate corollary of Lemma \ref{l1.2} is the following.

\begin{lemma}
\label{2.11}Suppose $\xi $ is Kuranishi data for a transversely holomorphic
trivialization of $M/\Delta $. Then for any 
\begin{equation*}
\alpha =\sum\nolimits_{\left| J\right| >0}\alpha _{J}t^{J},\quad \alpha
_{J}\in A^{0}\left( T_{1,0}\left( M_{0}\right) \right) ,
\end{equation*}
the family of Kuranishi data 
\begin{equation*}
\xi _{s}=\exp \left( \left[ s\alpha ,\ \right] \right) \left( \xi \right) -%
\frac{1-\exp \left( \left[ s\alpha ,\ \right] \right) }{\left[ s\alpha ,\
\right] }\left( s\left[ \bar{\partial },\alpha \right] \right) .
\end{equation*}
corresponds to a (formal) family of transversely holomorphic trivializations
of the same deformation $M/\Delta .$ If $\alpha $ is sufficiently small with
respect to a metric induced by a hermitian metric on $T_{1,0}\left(
M_{0}\right) $, then $\xi _{s}$ is convergent, i.e. corresponds to an actual
transversely holomorphic trivialization of $M/\Delta $, for sufficiently
small $s$. Also any transversely holomorphic trivialization ``sufficiently
close'' to $\xi $ can be obtained in this way.
\end{lemma}

We call $\xi $ \textit{integrable} if it comes from a transversely
holomorphic trrivialization of a deformation $M/\Delta $. Integrability is
given by convergence and the Newlander-Nirenberg condition (see \S \ref{A3}
of the Appendix). The proof of the existence of transversely holomorphic
trivializations shows that the space of them is path-connected. So we have:

\begin{lemma}
Two integrable series $\xi ^{\prime },\xi ^{\prime \prime }$ give
holomorphically equivalent deformations of $M_{0}$ over $\Delta $ if and
only if they can be connected by a path of transversely holomorphic
trivializations $\xi \left( s\right) $ for which there is a path of vector
fields 
\begin{equation*}
\alpha \left( s\right) =\sum\nolimits_{\left| J\right| >0}\alpha _{J}\left(
s\right) t^{J},\quad \alpha _{I}\left( s\right) \in A^{0,0}\left(
T_{1,0}\left( M_{0}\right) \right) 
\end{equation*}
such that 
\begin{equation*}
\frac{d\xi \left( s\right) }{du}=\bar{\partial }\alpha \left( s\right)
+\left[ \alpha \left( s\right) ,\xi \left( s\right) \right] .
\end{equation*}
\end{lemma}

We next need to understand the contents of Lemma \ref{2.5} in a more formal
way.

\section{The fundamental differential graded Lie algebra}

The ring of $C^{\infty }$-functions on a formal neighborhood of $M_{0}$ in $%
M $ is simply the set of formal power series 
\begin{equation*}
\sum\nolimits_{I,J}f_{I}t^{I}\overline{t^{J}}=\sum\nolimits_{I,J}\left(
f_{I}\circ \sigma \right) t^{I}\overline{t^{J}}
\end{equation*}
with each $f_{I}$ a $C^{\infty }$-function on $M_{0}$. Similarly via $\sigma 
$ the algebra $A^{*}\left( \left\langle M\right\rangle \right) $ of $%
C^{\infty }$-differentials on a formal neighborhood $\left\langle
M\right\rangle $ of $M_{0}$ in $M$ is identified with the algebra 
\begin{equation*}
\bigwedge \Bbb{C}\left[ dt,\overline{dt}\right] \otimes \left\{
\sum\nolimits_{I,J}\eta _{I,J}t^{I}\overline{t^{J}}\right\}
\end{equation*}
with each $\eta _{I,J}\in A^{*}\left( M_{0}\right) $ a $C^{\infty }$
-differential on $M_{0}$. So the algebra 
\begin{equation*}
A^{*}\left( \left\langle M\right\rangle /\left\langle \Delta \right\rangle
\right)
\end{equation*}
of relative differentials becomes 
\begin{equation}
\left\{ \sum\nolimits_{I,J}\eta _{I}t^{I}\overline{t^{J}}\right\} .
\label{2.8.1}
\end{equation}
We next divide the algebra $\left( \ref{2.8.1}\right) $ by the ideal 
\begin{equation*}
\overline{\ t}\cdot A^{*}\left( \left\langle M\right\rangle /\left\langle
\Delta \right\rangle \right)
\end{equation*}
to obtain a complex 
\begin{equation}
B^{*}\left( \left\langle M\right\rangle /\left\langle \Delta \right\rangle
\right) ,\overline{\partial }_{M}.  \label{2.8.2}
\end{equation}

The natural maps 
\begin{equation*}
\left( \sigma _{t}^{-1}\right) ^{*}:A^{*}\left( M_{t}\right) \rightarrow
A^{*}\left( M_{0}\right)
\end{equation*}
induce maps 
\begin{equation*}
A^{p,q}\left( M/\Delta \right) \overset{\left( \sigma _{t}^{-1}\right) ^{*}}{
\rightarrow }\Bbb{C}\left[ \left[ t,\overline{t}\right] \right] \otimes
\sum\nolimits_{p^{\prime }+q^{\prime }=p+q}A^{p^{\prime },q^{\prime }}\left(
M_{0}\right)
\end{equation*}
and so maps: 
\begin{equation*}
\begin{array}{ccc}
B^{p,q}\left( \left\langle M\right\rangle /\left\langle \Delta \right\rangle
\right) & \overset{\varphi ^{p,q}}{\rightarrow } & \frac{\Bbb{C}\left[
\left[ t,\overline{t}\right] \right] \otimes A^{p+q}\left( M_{0}\right) }{%
\overline{t}\Bbb{C}\left[ \left[ t,\overline{t}\right] \right] \otimes
A^{p+q}\left( M_{0}\right) } \\ 
&  & \updownarrow = \\ 
&  & \Bbb{C}\left[ \left[ t\right] \right] \otimes A^{p+q}\left( M_{0}\right)
\end{array}
\end{equation*}
In fact, by $\left( \ref{2.5.3}\right) $ we have injections 
\begin{equation}
\varphi ^{p,q}:B^{p,q}\left( \left\langle M\right\rangle /\left\langle
\Delta \right\rangle \right) \rightarrow \Bbb{C}\left[ \left[ t\right]
\right] \otimes \left( A^{p,q}\left( M_{0}\right) +A^{p-1,q+1}\left(
M_{0}\right) +\ldots \right)  \label{2.8.3}
\end{equation}
fitting together to give the isomorphism of $d$-exterior algebras 
\begin{equation}
\varphi :B^{*}\left( \left\langle M\right\rangle /\left\langle \Delta
\right\rangle \right) \rightarrow \Bbb{C}\left[ \left[ t\right] \right]
\otimes A^{*}\left( M_{0}\right) .  \label{2.8.4}
\end{equation}

\section{Properties of the operator $\overline{\partial }-L_{\xi }\label{3}$}

For 
\begin{equation*}
\xi =\sum\nolimits_{J,\left| J\right| >0}\xi _{J}\cdot t^{J}
\end{equation*}
define the operator 
\begin{eqnarray*}
\bar{D}_{\xi } &:&A_{M_{0}\times \Delta /\Delta }^{p,q}\rightarrow
A_{M_{0}\times \Delta /\Delta }^{p,q+1}+A_{M_{0}\times \Delta /\Delta
}^{p-1,q+2} \\
\omega &=&\sum\nolimits_{I}t^{I}\omega _{I}\longmapsto \left( \overline{
\partial }-L_{\xi }\right) \left( \omega \right)
\end{eqnarray*}
by the formula: 
\begin{equation}
\bar{D}_{_{\xi }}\left( \omega \right) :=\left( \bar{\partial }%
_{M_{0}}-\sum\nolimits_{J,\left| J\right| >0}t^{J}L_{\xi _{J}}\right) \left(
\omega \right) =\sum\nolimits_{I}\bar{\partial }_{M_{0}}\omega
_{I}t^{I}-\sum\nolimits_{I,J,\left| J\right| >0}L_{\xi _{J}}\left( \omega
_{I}\right) t^{I+J}.  \label{3.3.1}
\end{equation}
We have from $\left( \ref{2.7.5}\right) $ that 
\begin{equation}
\bar{D}_{_{\xi }}\left( d\omega \right) =-d\bar{D}_{_{\xi }}\left( \omega
\right) .  \label{3.3.2}
\end{equation}
From $\left( \ref{2.7.3}\right) $ it follows that 
\begin{equation}
\bar{D}_{_{\xi }}\left( \omega \wedge \eta \right) =\bar{D}_{_{\xi }}\left(
\omega \right) \wedge \eta +\left( -1\right) ^{\deg \omega }\omega \wedge 
\bar{D}_{_{\xi }}\left( \eta \right) .  \label{3.3.3}
\end{equation}

\begin{lemma}
\label{3.1}For the isomorphisms of $d$-exterior algebras 
\begin{equation*}
\varphi :B^{*}\left( \left\langle M\right\rangle /\left\langle \Delta %
\right\rangle \right) \rightarrow \Bbb{C}\left[ \left[ t\right] \right]
\otimes A^{*}\left( M_{0}\right) 
\end{equation*}
in $\left( \ref{2.8.4}\right) $ we have the identity 
\begin{equation*}
\varphi \circ \overline{\partial }_{M}=\left( \overline{\partial }%
_{M_{0}}-L_{\xi }\right) \circ \varphi .
\end{equation*}
\end{lemma}

\begin{proof}
We have from $\left( \ref{2.7.5}\right) $ that 
\begin{equation}
L_{\xi }\left( d\omega \right) =-d\left( L_{\xi }\left( \omega \right)
\right) .  \label{Llemma}
\end{equation}
The assertion in the lemma is local so we can work locally on a coordinate
disk $W$ in $M$ on which we have holomorphic coordinates 
\begin{equation*}
\left( \left( u^{l}\right) ,\left( t^{k}\right) \right) .
\end{equation*}
Also, by $\left( \ref{1.7}\right) $ and Lemma \ref{2.5}, the assertion holds
for functions 
\begin{equation*}
\varphi \left( \overline{\partial }_{M}f\right) =\left( \overline{\partial }
_{M_{0}}-L_{\xi }\right) \left( \varphi \left( f\right) \right) .
\end{equation*}

Then 
\begin{equation*}
\overline{\partial }_{M}\left( u^{l}\right) =0
\end{equation*}
so by $\left( \ref{Llemma}\right) $%
\begin{equation*}
\left( \overline{\partial }_{M_{0}}-L_{\xi }\right) \varphi \left(
du^{l}\right) =0.
\end{equation*}
On the other hand 
\begin{equation*}
\varphi \left( \overline{u^{l}}\right)
\end{equation*}
is a power series in $\overline{t}$ and so by $\left( \ref{Llemma}\right) $%
\begin{equation*}
\left( \overline{\partial }_{M_{0}}-L_{\xi }\right) \varphi \left( \overline{
du^{l}}\right) =0.
\end{equation*}
Since every element of $B^{*}\left( \left\langle M\right\rangle
/\left\langle \Delta \right\rangle \right) $ can be written locally as a sum
of terms given by a function times a wedge of $du^{l}$'s and $\overline{
du^{l}}$'s and both $\overline{\partial }_{M}$ and $\left( \overline{
\partial }_{M_{0}}-L_{\xi }\right) $ satisfy the Leibniz rule by $\left( \ref
{3.3.3}\right) ,$ the proof is complete.
\end{proof}

\begin{lemma}
\label{3.4}The ``$\left( 0,1\right) $'' tangent distribution given by the
image of 
\begin{eqnarray*}
\iota  &:&T_{0,1}\left( M_{0}\right) \rightarrow T\left( M_{0}\right)
\otimes \Bbb{C}\left[ \left[ t\right] \right]  \\
&&id.-\left\langle \left. \xi \right| id.\right\rangle 
\end{eqnarray*}
gives, via complex conjugation, an almost complex structure on $M_{0}\times 
\Delta $.

i) This almost complex structure is integrable, that is, come from a
(formal) deformation/trivialization of $M_{0}$ as in $\left( \ref{2.1}%
\right) $ and Definition \ref{2.2}, if and only if 
\begin{equation*}
\left[ \bar{D}_{_{\xi }}+\overline{\partial }_{M_{0}},\left\langle \left.
\xi \right| \ \right\rangle \right] =0.
\end{equation*}
ii) This almost complex structure is integrable if and only if 
\begin{equation*}
\bar{D}_{_{\xi }}\circ \bar{D}_{_{\xi }}=0.
\end{equation*}
\end{lemma}

\begin{proof}
i) 
\begin{equation*}
\bar{D}_{_{\xi }}=\overline{\partial }_{M_{0}}-L_{\xi }
\end{equation*}
and, by $\left( \ref{2.7.7}\right) ,$%
\begin{equation*}
\left[ L_{\xi },\left\langle \xi \right. \left| \ \right\rangle \right]
=\left\langle \left. \left[ \xi ,\xi \right] \right| \ \right\rangle
\end{equation*}
so that 
\begin{eqnarray*}
\left[ \bar{D}_{_{\xi }}+\overline{\partial }_{M_{0}},\left\langle \left.
\xi \right| \ \right\rangle \right] &=&\left[ 2\overline{\partial }
_{M_{0}}-L_{\xi },\left\langle \left. \xi \right| \ \right\rangle \right] \\
&=&\left( 2\left\langle \left. \overline{\partial }_{M_{0}}\xi \right| \
\right\rangle -\left\langle \left. \left[ \xi ,\xi \right] \right| \
\right\rangle \right) \\
&=&2\left\langle \left. \left( \overline{\partial }_{M_{0}}\xi -\frac{\left[
\xi ,\xi \right] }{2}\right) \right| \ \right\rangle .
\end{eqnarray*}
So by the classical integrability criterion in \S \ref{A3} the proof of i)
is complete.

ii) One direction is immediate from Lemma \ref{3.1} and the fact that 
\begin{equation*}
\overline{\partial }_{M}\circ \overline{\partial }_{M}=0.
\end{equation*}
For the other direction we compute (using the Einstein summation convention
and $\left( \ref{48.5}\right) $, $\left( \ref{2.7.2}\right) $, and $\left( 
\ref{2.7.6}\right) $): 
\begin{equation*}
\begin{tabular}{l}
$\bar{D}_{_{\xi }}\bar{D}_{_{\xi }}\left( \omega \right) =\bar{D}_{_{\xi
}}\left( \bar{\partial }\omega -t^{J}\left( \left\langle \xi _{J}\right.
\left| d\omega \right\rangle -d\left\langle \xi _{J}\right. \left| \omega
\right\rangle \right) \right) $ \\ 
$=\bar{D}_{_{\xi }}\bar{\partial }\omega -\bar{D}_{_{\xi }}\left(
\left\langle \xi _{J}\right. \left| d\omega \right\rangle -d\left\langle \xi
_{J}\right. \left| \omega \right\rangle \right) t^{J}$ \\ 
$=-\left( \left\langle \xi _{J}\right. \left| d\bar{\partial }\omega
\right\rangle -d\left\langle \xi _{J}\right. \left| \bar{\partial }\omega
\right\rangle \right) t^{J}$ \\ 
$-\bar{\partial }\left( \left\langle \xi _{J}\right. \left| d\omega
\right\rangle -d\left\langle \xi _{J}\right. \left| \omega \right\rangle
\right) t^{J}+t^{J}t^{K}L_{\xi _{K}}L_{\xi _{J}}\left( \omega \right) $ \\ 
$=\left( \left\langle \xi _{J}\right. \left| \bar{\partial }d\omega
\right\rangle +d\left\langle \xi _{J}\right. \left| \bar{\partial }\omega
\right\rangle \right) t^{J}$ \\ 
$-\left( \bar{\partial }\left\langle \xi _{J}\right. \left| d\omega
\right\rangle +d\bar{\partial }\left\langle \xi _{J}\right. \left| \omega
\right\rangle \right) t^{J}+t^{J}t^{K}L_{\xi _{K}}L_{\xi _{J}}\left( \omega
\right) $ \\ 
$=\left( -\left\langle \bar{\partial }\xi _{J}\right. \left| d\omega
\right\rangle -d\left\langle \bar{\partial }\xi _{J}\right. \left| \omega
\right\rangle \right) t^{J}+t^{J}t^{K}L_{\xi _{K}}L_{\xi _{J}}\left( \omega
\right) $ \\ 
$=-\left( L_{\bar{\partial }\xi }-\frac{1}{2}L_{\left[ \xi ,\xi \right]
}\right) \left( \omega \right) .$%
\end{tabular}
\end{equation*}
So again by the classical integrability criterion in \S \ref{A3} the proof
is complete.
\end{proof}

\begin{corollary}
\label{3.6}Suppose we have a deformation $M/\Delta $ given by Kuranishi data 
$\xi $ such that 
\begin{equation*}
\overline{\partial }_{M_{0}}\xi =0.
\end{equation*}
Then by integrability 
\begin{equation*}
\left[ \xi ,\xi \right] =0
\end{equation*}
and for 
\begin{equation*}
\overline{D}_{\sigma }=\overline{D}_{\xi }
\end{equation*}
we have 
\begin{equation*}
\left[ \overline{D}_{\sigma },\left\langle \left. \xi \right| \
\right\rangle \right] =0.
\end{equation*}
\end{corollary}

\section{Gauss-Manin connection\label{newGM}}

From \S \ref{A3} the $\left( 0,1\right) $-distribution $T_{0,1}\left(
M/\Delta \right) $ is given on $M_{0}\times \Delta $ as 
\begin{equation*}
\left( id.-\left\langle \left. \ \right| \xi \right\rangle \right) \left(
T_{0,1}\left( M_{0}\right) \right)
\end{equation*}
and the holomorphic cotangent bundle $T^{1,0}\left( M/\Delta \right) $ is 
\begin{equation*}
id.+\left\langle \left. \xi \right| \ \right\rangle \left( T^{1,0}\left(
M_{0}\right) \right) =e^{\left\langle \left. \xi \right| \ \right\rangle
}\left( T^{1,0}\left( M_{0}\right) \right) .
\end{equation*}
Notice that 
\begin{equation*}
\left\langle \left. \xi \right| \ \right\rangle :T^{0,1}\left( M_{0}\right)
\rightarrow T^{0,1}\left( M_{0}\right)
\end{equation*}
is the zero map. Also 
\begin{equation*}
\left\langle \left. \xi \right| \ \right\rangle \left( \alpha \wedge \beta
\right) =\left\langle \left. \xi \right| \alpha \right\rangle \wedge \beta
+\alpha \wedge \left\langle \left. \xi \right| \beta \right\rangle .
\end{equation*}
So, by induction 
\begin{eqnarray*}
\left\langle \left. \xi \right| \ \right\rangle ^{n}\left( \alpha \wedge
\beta \right) &=&\sum\nolimits_{r+s=n}\left( 
\begin{array}{c}
n \\ 
r
\end{array}
\right) \left\langle \left. \xi \right| \ \right\rangle ^{r}\left( \alpha
\right) \wedge \left\langle \left. \xi \right| \ \right\rangle ^{s}\left(
\beta \right) \\
\frac{\left\langle \left. \xi \right| \ \right\rangle ^{n}}{n!}\left( \alpha
\wedge \beta \right) &=&\sum\nolimits_{r+s=n}\frac{\left\langle \left. \xi
\right| \ \right\rangle ^{r}}{r!}\left( \alpha \right) \wedge \frac{%
\left\langle \left. \xi \right| \ \right\rangle ^{s}}{s!}\left( \beta \right)
\\
e^{\left\langle \left. \xi \right| \ \right\rangle }\left( \alpha \wedge
\beta \right) &=&e^{\left\langle \left. \xi \right| \ \right\rangle }\left(
\alpha \right) \wedge e^{\left\langle \left. \xi \right| \ \right\rangle
}\left( \beta \right) .
\end{eqnarray*}
So by multiplicativity, $\varphi B^{p,q}\left( \left\langle M\right\rangle
/\left\langle \Delta \right\rangle \right) $ is given by sections of 
\begin{equation}
e^{\left\langle \left. \xi \right| \ \right\rangle }T^{p,q}\left(
M_{0}\right) .  \label{2.9.1}
\end{equation}
Any such section $\eta $ can be written in the form 
\begin{equation*}
e^{\left\langle \left. \xi \right| \ \right\rangle }\eta ^{p,q}
\end{equation*}
where 
\begin{equation}
\eta ^{p,q}\in A_{M_{0}}^{p,q}\otimes \Bbb{C}\left[ \left[ t\right] \right] .
\label{2.9.1.1}
\end{equation}
is simply the summand of $\eta $ of type $\left( p,q\right) .$

Now the Gauss-Manin connection gives 
\begin{equation*}
\begin{array}{c}
L_{\tau _{k}}=L_{\tau _{k}}^{1,0}+L_{\tau _{k}}^{0,1}:B^{p,q}\left(
\left\langle M\right\rangle /\left\langle \Delta \right\rangle \right)
\rightarrow B^{p,q}\left( \left\langle M\right\rangle /\left\langle \Delta
\right\rangle \right) +B^{p-1,q+1}\left( \left\langle M\right\rangle
/\left\langle \Delta \right\rangle \right) . \\ 
\omega \longmapsto d\left\langle \left. \tau _{k}\right| \omega
\right\rangle +\left\langle \left. \tau _{k}\right| d\omega \right\rangle
\end{array}
\end{equation*}
where 
\begin{equation*}
\tau _{k}
\end{equation*}
is the vector field on $M$ which, under the isomorphism $\left( \ref{2.8.4}
\right) ,$ corresponds to $\frac{\partial }{\partial t_{k}}$. So by $\left( 
\ref{2.9.1}\right) :$

\begin{lemma}
\label{nl5}Under the isomorphism $\left( \ref{2.8.3}\right) $ $B^{p,q}\left(
\left\langle M\right\rangle /\left\langle \Delta \right\rangle \right) $
corresponds to the space 
\begin{equation*}
e^{\left\langle \left. \xi \right| \ \right\rangle }\left( A^{p,q}\left(
M_{0}\right) \otimes \Bbb{C}\left[ \left[ t\right] \right] \right) 
\end{equation*}
and the Gauss-Manin connection 
\begin{equation*}
L_{\tau _{k}}^{1,0}+L_{\tau _{k}}^{0,1}:B^{p,q}\left( \left\langle
M\right\rangle /\left\langle \Delta \right\rangle \right) \rightarrow
B^{p,q}\left( \left\langle M\right\rangle /\left\langle \Delta \right\rangle
\right) \oplus B^{p-1,q+1}\left( \left\langle M\right\rangle /\left\langle 
\Delta \right\rangle \right) 
\end{equation*}
is given by the operator 
\begin{eqnarray}
e^{\left\langle \left. \xi \right| \ \right\rangle }\left( T^{p,q}\left(
M_{0}\right) \right)  &\rightarrow &e^{\left\langle \left. \xi \right| \
\right\rangle }\left( T^{p,q}\left( M_{0}\right) \right) \oplus
e^{\left\langle \left. \xi \right| \ \right\rangle }\left( T^{p-1,q+1}\left(
M_{0}\right) \right) .  \label{2.9.2} \\
e^{\left\langle \left. \xi \right| \ \right\rangle }\tilde{\omega} &\mapsto
&e^{\left\langle \left. \xi \right| \ \right\rangle }\left( \frac{\partial 
\tilde{\omega}}{\partial t_{k}}\right) +e^{\left\langle \left. \xi \right| \
\right\rangle }\left\langle \left. \frac{\partial \xi }{\partial t_{k}}%
\right| \tilde{\omega}\right\rangle 
\end{eqnarray}
\end{lemma}

Furthermore, by Lemma \ref{3.4}i) we have 
\begin{equation}
\left[ \left( \overline{D}_{\sigma }+\overline{\partial }_{M_{0}}\right)
,e^{\left\langle \left. \xi \right| \ \right\rangle }\right] =0.
\label{misform}
\end{equation}

\begin{lemma}
\label{nl6}Under the isomorphism 
\begin{equation*}
e^{\left\langle \left. \xi \right| \ \right\rangle }:A^{p,q}\left(
M_{0}\right) \otimes \Bbb{C}\left[ \left[ t\right] \right] \rightarrow
\varphi \left( B^{p,q}\left( \left\langle M\right\rangle /\left\langle 
\Delta \right\rangle \right) \right) 
\end{equation*}
in $\left( \ref{2.8.4}\right) ,$ the operator 
\begin{equation*}
\overline{D}_{\sigma }:\varphi \left( B^{p,q}\left( \left\langle
M\right\rangle /\left\langle \Delta \right\rangle \right) \right)
\rightarrow \varphi \left( B^{p,q+1}\left( \left\langle M\right\rangle
/\left\langle \Delta \right\rangle \right) \right) 
\end{equation*}
corresponds to the operator 
\begin{equation*}
\overline{\partial }_{M_{0}}-L_{\xi }^{1,0}:A^{p,q}\left( M_{0}\right)
\otimes \Bbb{C}\left[ \left[ t\right] \right] \rightarrow A^{p,q+1}\left(
M_{0}\right) \otimes \Bbb{C}\left[ \left[ t\right] \right] .
\end{equation*}
\end{lemma}

\begin{proof}
Abbreviate 
\begin{equation*}
\overline{\partial }:=\overline{\partial }_{M_{0}}.
\end{equation*}
Notice that, by $\left( \ref{2.7.9}\right) $%
\begin{equation*}
\left[ L_{\xi }^{0,1},\left\langle \left. \xi \right| \ \right\rangle
\right] =0
\end{equation*}
so that 
\begin{equation}
\left\langle \left. \xi \right| \ \right\rangle \circ \left[ \overline{
\partial },\left\langle \left. \xi \right| \ \right\rangle \right] =\left[ 
\overline{\partial },\left\langle \left. \xi \right| \ \right\rangle \right]
\circ \left\langle \left. \xi \right| \ \right\rangle .  \label{comm}
\end{equation}

Thus 
\begin{equation*}
\left[ \overline{\partial },e^{\left\langle \left. \xi \right| \
\right\rangle }\right] =e^{\left\langle \left. \xi \right| \ \right\rangle
}\circ \left[ \overline{\partial },\left\langle \left. \xi \right| \
\right\rangle \right]
\end{equation*}
and we can compute 
\begin{eqnarray*}
e^{-\left\langle \left. \xi \right| \ \right\rangle }\circ \overline{
\partial }\circ e^{\left\langle \left. \xi \right| \ \right\rangle }
&=&e^{-\left\langle \left. \xi \right| \ \right\rangle }\circ \left[ 
\overline{\partial },e^{\left\langle \left. \xi \right| \ \right\rangle
}\right] +\overline{\partial } \\
&=&\left[ \overline{\partial },\left\langle \left. \xi \right| \
\right\rangle \right] +\overline{\partial } \\
&=&\overline{\partial }-L_{\xi }^{0,1}.
\end{eqnarray*}
By Lemma \ref{3.4}i) 
\begin{equation*}
e^{-\left\langle \left. \xi \right| \ \right\rangle }\circ \left( 2\overline{
\partial }-L_{\xi }\right) \circ e^{\left\langle \left. \xi \right| \
\right\rangle }=2\overline{\partial }-L_{\xi }
\end{equation*}
but by the above 
\begin{equation*}
e^{-\left\langle \left. \xi \right| \ \right\rangle }\circ \left( 2\overline{
\partial }-L_{\xi }\right) \circ e^{\left\langle \left. \xi \right| \
\right\rangle }=\left( \overline{\partial }-L_{\xi }^{0,1}\right)
+e^{-\left\langle \left. \xi \right| \ \right\rangle }\circ \left( \overline{
\partial }-L_{\xi }\right) \circ e^{\left\langle \left. \xi \right| \
\right\rangle }
\end{equation*}
so that 
\begin{eqnarray*}
2\overline{\partial }-L_{\xi }^{1,0}-L_{\xi }^{0,1} &=&\overline{\partial }
-L_{\xi }^{0,1}+e^{-\left\langle \left. \xi \right| \ \right\rangle }\circ
\left( \overline{\partial }-L_{\xi }\right) \circ e^{\left\langle \left. \xi
\right| \ \right\rangle } \\
\overline{\partial }-L_{\xi }^{1,0} &=&e^{-\left\langle \left. \xi \right| \
\right\rangle }\circ \left( \overline{\partial }-L_{\xi }\right) \circ
e^{\left\langle \left. \xi \right| \ \right\rangle }.
\end{eqnarray*}
\end{proof}

\begin{corollary}
\label{nl7}Suppose that $M/\Delta $ is a family of compact K\"{a}hler
manifolds. If 
\begin{equation*}
\eta \in \varphi \left( B^{p,q}\left( \left\langle M\right\rangle
/\left\langle \Delta \right\rangle \right) \right) 
\end{equation*}
is $d$-closed, then 
\begin{equation*}
\eta =e^{\left\langle \left. \xi \right| \ \right\rangle }\eta ^{p,q}
\end{equation*}
with 
\begin{equation*}
\partial \eta ^{p,q}=0.
\end{equation*}
\end{corollary}

\begin{proof}
Since $\left[ \varphi ,d\right] =0,$%
\begin{equation*}
\left( \partial +\overline{\partial }\right) \left( e^{\left\langle \left.
\xi \right| \ \right\rangle }\left( \eta ^{p,q}\right) \right) =0,
\end{equation*}
so that by type 
\begin{equation*}
\partial \eta ^{p,q}=0.
\end{equation*}
\end{proof}

\section{Obstruction class\label{3.5}}

Next, suppose that 
\begin{equation}
M_{\frak{A}}/\Delta _{\frak{A}}  \label{3.5.1}
\end{equation}
is an infinitesimal deformation of compact K\"{a}hler manifolds of
(relative) dimension $m$. From the formula 
\begin{equation*}
\bar{\partial }\left[ \xi ,\xi ^{\prime }\right] =\left[ \bar{\partial }\xi
,\xi ^{\prime }\right] +\left( -1\right) ^{\deg \xi }\left[ \xi ,\bar{%
\partial }\xi ^{\prime }\right]
\end{equation*}
in $\left( \ref{2.7.4}\right) $ and the integrability condition 
\begin{equation*}
\bar{\partial }\xi -\frac{1}{2}\left[ \xi ,\xi \right] \in A^{0,2}\left(
T_{M_{0}}\right) \otimes \frak{A}
\end{equation*}
given in \S \ref{A3}, we compute that, modulo $\frak{mA}$, 
\begin{eqnarray*}
2\bar{\partial }\left[ \xi ,\xi \right] &=&2\left[ \bar{\partial }\xi ,\xi
\right] -2\left[ \xi ,\bar{\partial }\xi \right] \\
&\equiv &\left[ \left[ \xi ,\xi \right] ,\xi \right] -\left[ \xi ,\left[ \xi
,\xi \right] \right] =0
\end{eqnarray*}
so that 
\begin{equation}
\left[ \xi ,\xi \right] \in H^{2}\left( T_{M_{0}}\right) \otimes \frac{\frak{%
A}}{\frak{mA}}.  \label{3.5.2}
\end{equation}
To see that this is the obstruction to extending $M_{\frak{A}}/\Delta _{%
\frak{A}}$ to a family 
\begin{equation*}
M_{\frak{A}^{\prime }}/\Delta _{\frak{A}^{\prime }}
\end{equation*}
for some 
\begin{equation*}
\frak{A}\supseteq \frak{A}^{\prime }\supseteq \frak{mA}
\end{equation*}
suppose that 
\begin{equation*}
\bar{\partial }\alpha -\left( \bar{\partial }\xi -\frac{1}{2}\left[ \xi ,\xi
\right] \right) \in A^{0,2}\left( T_{M_{0}}\right) \otimes \frak{A}^{\prime }
\end{equation*}
for some 
\begin{equation*}
\alpha \in A^{0,1}\left( T_{M_{0}}\right) \otimes \frak{A.}
\end{equation*}
Then 
\begin{equation*}
\bar{\partial }\left( \xi -\alpha \right) -\left[ \left( \xi -\alpha \right)
,\left( \xi -\alpha \right) \right] \in A^{0,2}\left( T_{M_{0}}\right)
\otimes \frak{A}^{\prime }.
\end{equation*}

\part{Cohomology and obstructions\label{part4}}

\section{Annihilation of obstructions by cohomology classes}

\begin{theorem}
\label{4.1}Given a K\"{a}hler deformation 
\begin{equation*}
M_{\frak{A}}/\Delta _{\frak{A}}
\end{equation*}
of a compact K\"{a}hler manifold $M_{0}$ as in \ref{3.5} and obstruction
class $\left[ \xi ,\xi \right] $ as in $\left( \ref{3.5.2}\right) $ and
given any element 
\begin{equation*}
\omega _{0}\in H^{p,q}\left( M_{0}\right) ,
\end{equation*}
the value of the pairing 
\begin{equation*}
\left\langle \left. \left[ \xi ,\xi \right] \right| \omega _{0}\right\rangle
\in H^{p-1,q+2}\left( M_{0}\right) \otimes \frac{\frak{A}}{\frak{mA}}
\end{equation*}
is zero.
\end{theorem}

\begin{proof}
We realize the deformation $M_{\frak{A}}/\Delta _{\frak{A}}$ as Kuranishi
data $\xi $ on $M_{0}\times \Delta $ such that 
\begin{equation*}
\overline{\partial }_{M_{0}}\xi -\frac{1}{2}\left[ \xi ,\xi \right] \in
A^{0,2}\left( T_{M_{0}}\right) \otimes \frak{A}.
\end{equation*}
Abbreviate 
\begin{equation*}
\overline{\partial }:=\overline{\partial }_{M_{0}}
\end{equation*}
throughout the remainder of the proof. From \S \ref{3} and \S \ref{newGM}
there is an extension of $\omega _{0}$ to a form $\omega $ of type $\left(
p,q\right) $ on $M_{0}\times \Delta $ such that 
\begin{equation*}
\left. d\left( e^{\left\langle \left. \xi \right| \ \right\rangle }\omega
\right) \right| _{M_{0}\times \Delta _{\frak{A}}}=0
\end{equation*}
and 
\begin{equation*}
\left. \left( \overline{\partial }-L_{\xi }^{1,0}\right) \left( \omega
\right) \right| _{M_{0}\times \Delta _{\frak{A}}}=0.
\end{equation*}
Thus by type 
\begin{equation*}
\left. \partial \omega \right| _{M_{0}\times \Delta _{\frak{A}}}=0.
\end{equation*}

By $\left( \ref{2.7.9}\right) $ and the fact that $\omega $ is $\partial $
-closed, we have modulo $\frak{mA}$ that 
\begin{eqnarray*}
-\left\langle \left. \left[ \xi ,\xi \right] \right| \omega \right\rangle
&=&\left\langle \left. \xi \right| L_{\xi }^{1,0}\omega \right\rangle
-L_{\xi }^{1,0}\left\langle \left. \xi \right| \omega \right\rangle \\
&=&-2\left\langle \left. \xi \right| \partial \left\langle \left. \xi
\right| \omega \right\rangle \right\rangle +\partial \left\langle \left. \xi
\right| \left\langle \left. \xi \right| \omega \right\rangle \right\rangle
\end{eqnarray*}
so that 
\begin{eqnarray}
\frac{1}{2}\left\langle \left. \left[ \xi ,\xi \right] \right| \omega
\right\rangle &=&\left\langle \left. \xi \right| \partial \left\langle
\left. \xi \right| \omega \right\rangle \right\rangle -\frac{1}{2}\partial
\left\langle \left. \xi \right| \left\langle \left. \xi \right| \omega
\right\rangle \right\rangle  \label{4.1.2} \\
&=&-\left\langle \left. \xi \right| L_{\xi }^{1,0}\omega \right\rangle -%
\frac{1}{2}\partial \left\langle \left. \xi \right| \left\langle \left. \xi
\right| \omega \right\rangle \right\rangle .  \notag
\end{eqnarray}
Since 
\begin{equation*}
\left( \overline{\partial }-L_{\xi }^{1,0}\right) \omega \in A^{*}\left(
M_{0}\right) \otimes \frak{A},
\end{equation*}
and 
\begin{equation*}
\xi \in A^{0,1}\left( T_{M_{0}}\right) \otimes \frak{m},
\end{equation*}
we can use $\left( \ref{4.1.2}\right) $ to compute modulo $\frak{mA}$: 
\begin{eqnarray*}
\overline{\partial }\left\langle \left. \xi \right| \omega \right\rangle
&=&\left\langle \left. \overline{\partial }\xi \right| \omega \right\rangle
+\left\langle \left. \xi \right| \overline{\partial }\omega \right\rangle \\
&\equiv &\left\langle \left. \overline{\partial }\xi \right| \omega
\right\rangle +\left\langle \left. \xi \right| L_{\xi }^{1,0}\omega
\right\rangle \\
&\equiv &\left\langle \left. \overline{\partial }\xi \right| \omega
\right\rangle -\frac{1}{2}\left\langle \left. \left[ \xi ,\xi \right]
\right| \omega \right\rangle -\frac{1}{2}\partial \left\langle \left. \xi
\right| \left\langle \left. \xi \right| \omega \right\rangle \right\rangle \\
&\equiv &\left\langle \left. \overline{\partial }\xi -\frac{1}{2}\left[ \xi
,\xi \right] \right| \omega \right\rangle -\frac{1}{2}\partial \left\langle
\left. \xi \right| \left\langle \left. \xi \right| \omega \right\rangle
\right\rangle \\
&\equiv &\left\langle \left. \overline{\partial }\xi -\frac{1}{2}\left[ \xi
,\xi \right] \right| \omega _{0}\right\rangle -\frac{1}{2}\partial
\left\langle \left. \xi \right| \left\langle \left. \xi \right| \omega
\right\rangle \right\rangle .
\end{eqnarray*}
That is 
\begin{equation*}
\overline{\partial }\left\langle \left. \xi \right| \omega \right\rangle +%
\frac{1}{2}\partial \left\langle \left. \xi \right| \left\langle \left. \xi
\right| \omega _{0}\right\rangle \right\rangle \equiv \overline{\partial }%
\left\langle \left. \xi \right| \omega _{0}\right\rangle -\frac{1}{2}%
\left\langle \left. \left[ \xi ,\xi \right] \right| \omega _{0}\right\rangle
\end{equation*}
and so modulo $\frak{mA}$ we have that 
\begin{equation*}
\left\langle \left. \left[ \xi ,\xi \right] \right| \omega _{0}\right\rangle
\equiv 0\in H^{p-1,q+2}\left( M_{0}\right) \otimes \frac{\Bbb{C}\left[
\left[ t\right] \right] }{\frak{mA}}.
\end{equation*}
\end{proof}

\begin{corollary}
Let $M_{0}$ be a compact K\"{a}hler manifold such that 
\begin{equation*}
\omega _{M_{0}}=\mathcal{O}_{M_{0}}.
\end{equation*}
Then $M_{0}$ is unobstructed, that is, the subspace 
\begin{equation*}
Obs\subseteq H^{2}\left( T_{M_{0}}\right) 
\end{equation*}
generated by elements occurring in $\left( \ref{3.5.2}\right) $, is zero.
\end{corollary}

\begin{proof}
The map 
\begin{equation*}
H^{2}\left( T_{M_{0}}\right) \rightarrow Hom\left( H^{0}\left( \omega
_{M_{0}}\right) ,H^{2}\left( \Omega _{M_{0}}^{n-1}\right) \right)
\end{equation*}
is an isomorphism in this case.
\end{proof}

\section{Deformations of submanifolds\label{4.2}}

Suppose again that we have a deformation $M/\Delta $ of the compact
K\"{a}hler manifold $M_{0}$ as in $\left( \ref{2.1}\right) $. Suppose
further that we are given a complex submanifold 
\begin{equation}
Y_{0}\subseteq M_{0}.  \label{4.2.1}
\end{equation}
If $Y_{0}$ is not compact we require that it has regular boundary, that is,
there is a system of open submanifolds 
\begin{equation*}
Y_{0}\subseteq Y_{0}\left[ \varepsilon \right) \subseteq M_{0},\varepsilon
\in \left( 0,1\right)
\end{equation*}
such that 
\begin{equation*}
\varepsilon \leq \varepsilon ^{\prime }\Longrightarrow Y_{0}\left[
\varepsilon \right) \subseteq Y_{0}\left[ \varepsilon ^{\prime }\right)
\end{equation*}
\begin{equation*}
\bigcap\nolimits_{\varepsilon }Y_{0}\left[ \varepsilon \right)
-Y_{0}=\partial Y_{0}
\end{equation*}
where $\partial Y_{0}$ is a compact real submanifold of real codimension one
lying inside each $Y_{0}\left[ \varepsilon \right) ,$ and 
\begin{equation*}
Y_{0}\left[ \varepsilon \right) -Y_{0}\underset{diffeo.}{\sim }\partial
Y_{0}\times \left[ 0,\varepsilon \right) .
\end{equation*}
We take cohomology on $Y_{0}$ to mean the direct limit of cohomology groups
on $Y_{0}\left[ \varepsilon \right) .$

Suppose now that we are given a deformation 
\begin{equation}
Y_{\frak{A}}/\Delta _{\frak{A}}\subseteq M_{\frak{A}}/\Delta _{\frak{A}}
\label{4.2.2}
\end{equation}
of $Y_{0}$ for some ideal $\frak{A}\subseteq \frak{m\subseteq }\Bbb{C}\left[
\left[ t\right] \right] $. (In the non-compact case, we take this to mean a
deformation of $Y_{0}\left[ \varepsilon \right) $ over $\Delta _{\frak{A}}$
for some $\varepsilon >0$.) We proceed exactly as in Lemma \ref{2.3} to
construct a transversely holomorphic trivialization of the deformation $%
\left( \ref{2.1}\right) $ except that we now additionally require that, in
each of our initial choices of holomorphic coordinates $u_{W}$, $Y_{\frak{A}%
} $ is (locally) defined in $M_{\frak{A}}$ by setting a subset of the
coordinates $u_{W}$ equal to zero. One then obtains that 
\begin{equation}
Y_{\frak{A}}\subseteq \sigma ^{-1}\left( Y_{0}\right)  \label{4.2.3}
\end{equation}
so that the restriction of the Kuranishi datum on $M_{\frak{A}}$ is the
Kuranishi datum for $Y_{\frak{A}}$, that is, 
\begin{equation}
\left. \xi \right| _{Y_{0}\times \Delta }=\xi ^{\prime }+\xi ^{\prime \prime
}\in \left( A^{0,1}\left( T_{Y_{0}}\right) \otimes \frak{m}+A^{0,1}\left(
\left. T_{M_{0}}\right| _{Y_{0}}\right) \otimes \frak{A}\right) .
\label{4.2.4}
\end{equation}
Thus 
\begin{equation*}
\left. \left[ \xi ,\xi \right] \right| _{Y_{0}\times \Delta }\in \left(
A^{0,2}\left( T_{Y_{0}}\right) \otimes \frak{m}^{2}+A^{0,2}\left( \left.
T_{M_{0}}\right| _{Y_{0}}\right) \otimes \frak{mA}\right)
\end{equation*}
and so, if $N_{*\backslash *}$ denotes the normal bundle, 
\begin{equation*}
\left. \left[ \xi ,\xi \right] \right| _{Y_{0}\times \Delta }\in
A^{0,2}\left( N_{Y_{0}\backslash M_{0}}\right) \otimes \frak{mA}
\end{equation*}
and the integrability condition (see \S \ref{A3}) for the deformation of $%
M_{0}$ then implies that 
\begin{equation*}
\left. \overline{\partial }\xi \right| _{Y_{0}\times \Delta }=\frac{1}{2}%
\left. \left[ \xi ,\xi \right] \right| _{Y_{0}\times \Delta }\in
A^{0,2}\left( N_{Y_{0}\backslash M_{0}}\right) \otimes \frak{mA}
\end{equation*}
so that 
\begin{equation}
\left\{ \xi ^{\prime \prime }\right\} \in H^{1}\left( N_{Y_{0}\backslash
M_{0}}\right) \otimes \frac{\frak{A}}{\frak{mA}}  \label{4.2.5}
\end{equation}
is the obstruction to extending the deformation $\left( \ref{4.2.2}\right) $
to a family 
\begin{equation*}
Y_{\frak{mA}}/\Delta _{\frak{mA}}\subseteq M_{\frak{mA}}/\Delta _{\frak{mA}}.
\end{equation*}
Indeed, suppose the element $\left( \ref{4.2.5}\right) $ vanishes modulo $%
\frak{A}^{\prime }$ for some $\frak{A}\supseteq \frak{A}^{\prime }\supseteq 
\frak{mA}$, that is, 
\begin{equation*}
\overline{\partial }\alpha \equiv -\xi \in A^{0,1}\left( N_{Y_{0}\backslash
M_{0}}\right) \otimes \frac{\frak{A}^{\prime }}{\frak{mA}}
\end{equation*}
for some 
\begin{equation*}
\alpha \in A^{0}\left( T_{M_{0}}\right) \otimes \frak{A}.
\end{equation*}
Writing the formula 
\begin{eqnarray*}
\tilde{\xi} &=&\exp \left( \left[ \alpha ,\ \right] \right) \left( \xi
\right) -\frac{1-\exp \left( \left[ \alpha ,\ \right] \right) }{\left[
\alpha ,\ \right] }\left( \bar{\partial }\alpha \right) \\
&=&\xi +\left[ \alpha ,\xi \right] +\bar{\partial }\alpha +\ldots ,
\end{eqnarray*}
we have by Lemma \ref{l1.2} that

i) $\tilde{\xi}$ is associated to some trivialization $F_{\tilde{\sigma}
}=\left( \tilde{\sigma},\pi \right) $ of the same deformation $\left( \ref
{2.1}\right) $,

\noindent and

ii) 
\begin{equation*}
\left. \tilde{\xi}\right| _{Y_{0}\times \Delta }\in A^{0,1}\left(
T_{Y_{0}}\right) \otimes \frak{m}+A^{0,1}\left( \left. T_{M_{0}}\right|
_{Y_{0}}\right) \otimes \frak{A}^{\prime }
\end{equation*}
so that $\tilde{\sigma}^{-1}\left( Y_{0}\times \Delta _{\frak{A}^{\prime
}}\right) $ is a complex submanifold of $M_{\frak{A}^{\prime }}/\Delta _{%
\frak{A}^{\prime }}$.

\section{Semiregularity}

\begin{theorem}
\label{4.3}Given a deformation $M/\Delta $ of a compact K\"{a}hler manifold $%
M_{0}$, and a deformation $Y_{\frak{A}}$ of $Y_{0}\subseteq M_{0}$ over $%
\Delta _{\frak{A}}$ as in \ref{4.2} with $Y_{0}$ compact. Suppose further
that the sub-Hodge-structure 
\begin{equation*}
K_{0}^{r}=\sum K_{0}^{p,q}=\ker \left( H^{r}\left( M_{0}\right) \rightarrow
H^{r}\left( Y_{0}\right) \right) 
\end{equation*}
deforms over $\Delta $ for some $r$. Then, for 
\begin{equation*}
\omega _{0}\in K_{0}^{p+1,q-1}
\end{equation*}
and for the representative $\xi ^{\prime \prime }$ of obstruction class in $%
\left( \ref{4.2.5}\right) $, 
\begin{equation*}
\left. \left\langle \left. \xi ^{\prime \prime }\right| \omega
_{0}\right\rangle \right| _{Y_{0}}=0\in H^{p,q}\left( Y_{0}\right) \otimes 
\frac{\frak{A}}{\frak{mA}}.
\end{equation*}
\end{theorem}

\begin{proof}
We use $\sigma ^{-1}\left( Y_{0}\right) $ as the representative of the
prolongation of $Y_{0}$ via the Gauss-Manin connection. Referring to Lemma 
\ref{nl5} we prolong $\omega _{0}$ to a class 
\begin{equation*}
\omega \in A^{p,q}\left( M_{0}\right) \otimes \Bbb{C}\left[ \left[ t\right]
\right]
\end{equation*}
which is $\partial $- and $\left( \overline{\partial }-L_{\xi }^{1,0}\right) 
$-closed and is such that 
\begin{equation*}
e^{\left\langle \left. \xi \right| \ \right\rangle }\left( \omega \right)
\end{equation*}
is $d$ and $\overline{D}_{\sigma }$-closed. By hypothesis 
\begin{equation*}
\left. e^{\left\langle \left. \xi \right| \ \right\rangle }\left( \omega
\right) \right| _{Y_{0}\times \Delta _{\frak{A}}}
\end{equation*}
is cohomologous to zero. Since 
\begin{equation*}
\left. \varphi ^{-1}\omega \right| _{Y_{\frak{A}}}
\end{equation*}
is cohomologous to zero, by the $\partial \overline{\partial }$-lemma
(Corollary 5.4 of \cite{D}$)$ there is 
\begin{equation*}
\gamma \in B^{p,q-2}\left( \left\langle M\right\rangle /\left\langle \Delta
\right\rangle \right)
\end{equation*}
such that 
\begin{equation*}
\left. \partial _{M}\overline{\partial }_{M}\gamma \right| _{Y_{\frak{A}%
}}=\left. \varphi ^{-1}\omega \right| _{Y_{\frak{A}}}.
\end{equation*}
Adjusting the choice of $\varphi ^{-1}\omega $ by $\partial _{M}\overline{%
\partial }_{M}\gamma $ we can assume that 
\begin{equation*}
\left. \varphi ^{-1}\omega \right| _{Y_{\frak{A}}}=0
\end{equation*}
as a cocycle. Thus 
\begin{equation}
\left. e^{\left\langle \left. \xi \right| \ \right\rangle }\left( \omega
\right) \right| _{Y_{0}\times \Delta }\in A^{p+q}\left( Y_{0}\right) \otimes 
\frak{A}  \notag
\end{equation}
and 
\begin{equation*}
\left\{ \left. e^{\left\langle \left. \xi \right| \ \right\rangle }\left(
\omega \right) \right| _{Y_{0}\times \Delta }\right\} =0\in H^{p+q}\left(
Y_{0}\right) \otimes \Bbb{C}\left[ \left[ t\right] \right] .
\end{equation*}

On the other hand, by $\left( \ref{4.2.4}\right) $%
\begin{equation*}
\left. \xi \right| _{Y_{0}\times \Delta }=\xi ^{\prime }+\xi ^{\prime \prime
}\in \left( A^{0,1}\left( T_{Y_{0}}\right) \otimes \frak{m}+A^{0,1}\left(
\left. T_{M_{0}}\right| _{Y_{0}}\right) \otimes \frak{A}\right) .
\end{equation*}
So modulo $\frak{mA}$%
\begin{equation*}
\left. \left\langle \left. \xi \right| \omega \right\rangle \right|
_{Y_{0}\times \Delta }=\left. \left\langle \left. \xi ^{\prime \prime
}\right| \omega _{0}\right\rangle \right| _{Y_{0}\times \Delta }
\end{equation*}
is the obstruction class to extending $Y_{\frak{A}}$ to a subscheme $Y_{%
\frak{mA}}$.

Now writing 
\begin{equation*}
e^{\left\langle \left. \xi \right| \ \right\rangle }\left( \omega \right)
=e^{\left\langle \left. \xi ^{\prime }+\xi ^{\prime \prime }\right| \
\right\rangle }\left( \omega _{0}+t\cdot \ldots \right)
\end{equation*}
and restricting to $Y_{0}\times \Delta $ and working modulo $\frak{mA}$, we
have 
\begin{equation*}
\left. e^{\left\langle \left. \xi \right| \ \right\rangle }\left( \omega
\right) \right| _{Y_{0}\times \Delta }=\omega +\left\langle \left. \xi
^{\prime \prime }\right| \omega _{0}\right\rangle .
\end{equation*}
Since $\left. e^{\left\langle \left. \xi \right| \ \right\rangle }\left(
\omega \right) \right| _{Y_{0}\times \Delta }$ is cohomologous to zero, its
projection $\left\{ \left\langle \left. \xi ^{\prime \prime }\right| \omega
_{0}\right\rangle \right\} $ into $H^{p,q}\left( Y_{0}\right) \otimes \frac{%
\frak{A}}{\frak{mA}}$ must be zero, and so the proof is complete.
\end{proof}

\section{Semiregularity, curvilinear version\label{App}}

In the case in which $Y_{0}$ is only relatively compact, that is, proper
over $Y_{0}^{\prime }$, we can refine Theorem \ref{4.3} somewhat, a fact
which will be useful in applications. However to achieve this strengthening,
we must restrict our attention to so-called \textit{curvilinear}
deformations:

\begin{theorem}
\label{4.4}Suppose $\dim \Delta =1$ and that $M/\Delta $ is a deformation of
the complex projective manifold $M_{0}$. Suppose 
\begin{equation*}
p^{\prime }:Y_{0}\rightarrow Y_{0}^{\prime }
\end{equation*}
is a proper family of submanifolds of $M_{0}$ of fiber dimension $q$ over a
smooth (not necessarily compact) base $Y_{0}^{\prime }$ of dimension $\geq p$
. Suppose further that the family $Y_{0}/Y_{0}^{\prime }$ deforms with $M_{0}
$ to a family $Y_{n}/Y_{n}^{\prime }$ over the Artinian scheme $\Delta %
_{n}\subseteq \Delta $ associated to the ideal $\left\{ t^{n+1}\right\} $
and that 
\begin{equation*}
\left\{ \tilde{\omega}\right\} \in H^{p+q+1,q-1}\left( M/\Delta \right) 
\end{equation*}
lies in the kernel of the composition 
\begin{equation*}
H^{p+q+1,q-1}\left( M/\Delta \right) \overset{L_{\tau }^{0,1}}{\rightarrow }%
H^{p+q,q}\left( M/\Delta \right) \overset{\left( pull-back\right) }{%
\longrightarrow }R^{q}p_{*}^{\prime }\left( \Omega _{Y_{n}/\Delta %
_{n}}^{p+q}\right) \rightarrow \Omega _{Y_{n}^{\prime }}^{p}
\end{equation*}
induced by the Gauss-Manin connection and integration over the fiber. Let 
\begin{equation*}
\xi ^{\prime \prime }\in R^{1}p_{*}^{\prime }\left( N_{Y_{0}\backslash
Y_{0}^{\prime }\times M_{0}}\right) 
\end{equation*}
be the obstruction class measuring extendability of $Y_{n}/Y_{n}^{\prime }$
to a family $Y_{n+1}/Y_{n+1}^{\prime }.$ Then, for $\omega _{0}=\left. 
\tilde{\omega}\right| _{M_{0}},$%
\begin{equation*}
\left. \left\langle \left. \xi ^{\prime \prime }\right| \omega
_{0}\right\rangle \right| _{Y_{0}}\in R^{q}p_{*}^{\prime }\left( \Omega %
_{Y_{0}}^{p+q}\right) \otimes \frac{_{\left\{ t^{n+1}\right\} }}{_{\left\{
t^{n+2}\right\} }}
\end{equation*}
goes to zero in 
\begin{equation*}
\Omega _{Y_{0}^{\prime }}^{p}\otimes \frac{_{\left\{ t^{n+1}\right\} }}{%
_{\left\{ t^{n+2}\right\} }}
\end{equation*}
under the map 
\begin{equation*}
R^{q}p_{*}^{\prime }\left( \Omega _{Y_{0}}^{p+q}\right) \rightarrow \Omega %
_{Y_{0}^{\prime }}^{p}
\end{equation*}
induced by the Leray spectral sequence. (This last map is commonly called 
\textit{integration over the fiber}.)
\end{theorem}

\begin{proof}
The statement is local on $Y_{0}^{\prime }$ so we can assume that it is an
analytic polydisk. Also we can assume that 
\begin{equation*}
\dim Y_{0}^{\prime }=p
\end{equation*}
since we can always restrict the original family to one of that dimension.
Thus the map 
\begin{equation*}
R^{q}p_{*}^{\prime }\left( \Omega _{Y_{n}/\Delta _{n}}^{p+q}\right)
\rightarrow \Omega _{Y_{n}^{\prime }}^{p}
\end{equation*}
can be assumed to be an isomorphism. We can further assume that $Y_{0}$ is
actually imbedded in $M_{0}$ since, if not, complete $Y_{0}^{\prime }$ as
follows. The varieties parametrized by $Y_{0}^{\prime }$ lie in some
component $\mathrm{Hilb}/\Delta $ of the relative Hilbert scheme of some
projective space (over $\Delta $) containing $M/\Delta $. Thus $%
Y_{0}^{\prime }$ is always an analytic open set in some projective variety $%
V_{0}^{\prime }$ such that the family $Y_{0}/Y_{0}^{\prime }$ is induced by
the restriction of a mapping 
\begin{equation*}
V_{0}^{\prime }\rightarrow \mathrm{Hilb}_{0}.
\end{equation*}
We replace $M/\Delta $ with 
\begin{equation*}
M\times _{\Delta }P/\Delta
\end{equation*}
where $P/\Delta $ is a projective space containing $\mathrm{Hilb}/\Delta $
and replace $Y_{0}$ by its image under the map 
\begin{equation*}
Y_{0}\rightarrow M_{0}\times Y_{0}^{\prime }\rightarrow M_{0}\times P_{0}.
\end{equation*}

Since we can assume that $Y_{0}$ is imbedded in $M_{0}$, we can construct a
trivialization so that $\sigma ^{-1}\left( Y_{0}\right) $ represents the
prolongation of $Y_{0}$ via the Gauss-Manin connection, in fact, in such a
way that the diagram 
\begin{equation}
\begin{array}{ccc}
Y_{n} & \overset{F_{\sigma }}{\longrightarrow } & Y_{0}\times \Delta _{n} \\ 
\downarrow p^{\prime } &  & \downarrow p_{0}^{\prime }\times 1_{\Delta _{n}}
\\ 
Y_{n}^{\prime } & \overset{G}{\longrightarrow } & Y_{0}^{\prime }\times
\Delta _{n}
\end{array}
\label{cd1}
\end{equation}
is commutative for some analytic isomorphism $G$.

Now 
\begin{equation*}
\left\{ \left. L_{\tau }^{0,1}\tilde{\omega}\right| _{Y_{n}}\right\} =0\in
R^{p+q}p_{*}^{\prime }\left( \Omega _{Y_{n}/Y_{n}^{\prime }}^{q}\right)
\cong \Omega _{Y_{0}^{\prime }\times \Delta _{n}}^{p}
\end{equation*}
and we let 
\begin{eqnarray*}
\varphi \left( \tilde{\omega}\right) &=&:\omega =e^{\left\langle \left. \xi
\right| \ \right\rangle }\left( \omega _{0}+t\cdot \ldots \right) \\
&\in &e^{\left\langle \left. \xi \right| \ \right\rangle }\left(
A^{p+1,q-1}\left( M_{0}\right) \otimes \Bbb{C}\left[ \left[ t\right] \right]
\right) .
\end{eqnarray*}
By Lemma \ref{nl5} and $\left( \ref{2.7.00}\right) $%
\begin{equation*}
\varphi L_{\tau }^{0,1}\left( \tilde{\omega}\right) =e^{\left\langle \left.
\xi \right| \ \right\rangle }\left\langle \left. \frac{\partial \xi }{
\partial t}\right| \omega \right\rangle
\end{equation*}
and so we have 
\begin{equation}
\left\{ \left. e^{\left\langle \left. \xi \right| \ \right\rangle
}\left\langle \left. \frac{\partial \xi }{\partial t}\right| \omega
\right\rangle \right| _{Y_{0}\times \Delta _{n}}\right\} =0\in \Omega
_{Y_{0}^{\prime }\times \Delta _{n}}^{p}.  \label{xxx}
\end{equation}
On the other hand, using $\left( \ref{cd1}\right) ,$%
\begin{equation*}
\left. \xi \right| _{Y_{0}\times \Delta }=\xi ^{\prime }+\xi ^{\prime \prime
}\in \left( A^{0,1}\left( T_{Y_{0}/Y_{0}^{\prime }}\right) \otimes \left\{
t\right\} +A^{0,1}\left( \left. T_{M_{0}}\right| _{Y_{0}}\right) \otimes
\left\{ t^{n+1}\right\} \right) .
\end{equation*}
where $T_{Y_{0}/Y_{0}^{\prime }}$ denotes the relative tangent space, that
is, the kernel of the projection 
\begin{equation*}
\left( p_{0}\right) _{*}:T_{Y_{0}}\rightarrow p_{0}^{*}T_{Y_{0}^{\prime }}.
\end{equation*}
But then 
\begin{equation*}
\left. \left\langle \left. \xi ^{\prime }\right| \omega \right\rangle
\right| _{Y_{0}\times \Delta }=0
\end{equation*}
since 
\begin{equation*}
\left. \omega \right| _{Y_{0}\times \Delta }=0
\end{equation*}
by dimension. So all summands in the expression for 
\begin{equation*}
\left. e^{\left\langle \left. \xi \right| \ \right\rangle }\left\langle
\left. \frac{\partial \xi }{\partial t}\right| \omega \right\rangle \right|
_{Y_{0}\times \Delta }
\end{equation*}
must be zero unless they involve $\xi ^{\prime \prime }$ or $\frac{\partial
\xi ^{\prime \prime }}{\partial t}$. So, writing 
\begin{equation*}
\xi ^{\prime \prime }=\xi _{n+1}^{\prime \prime }t^{n+1}+\ldots ,
\end{equation*}
the first possibly non-zero term is the coefficient of $t^{n}$ and that
coefficient is 
\begin{equation*}
\left\langle \left. \left( n+1\right) \xi _{n+1}^{\prime \prime }\right|
\omega _{0}\right\rangle .
\end{equation*}
Therefore by $\left( \ref{xxx}\right) $%
\begin{equation*}
\left\{ \left\langle \left. \left( n+1\right) \xi _{n+1}^{\prime \prime
}\right| \omega _{0}\right\rangle \right\} =0.
\end{equation*}
But $\xi _{n+1}^{\prime \prime }$ is the obstruction class for entending $%
Y_{n}/Y_{n}^{\prime }$ to a family over $Y_{0}^{\prime }\times \Delta _{n+1}$
.
\end{proof}

The main purpose of the sequel paper, ``\textit{Cohomology and Obstructions
II: Curves on Calabi-Yau threefolds''} will be to give some concrete
geometric applications of Theorem \ref{4.4}.

\section{Deformations of a pair\label{5}}

Finally we consider the case of a pair $\left( M_{0},Y_{0}\right) $ where $%
M_{0}$ is a complex manifold and $Y_{0}$ is a locally closed submanifold as
in \S \ref{4.2}. Suppose we have a deformation 
\begin{equation*}
Y_{\frak{A}}\backslash M_{\frak{A}}/\Delta _{\frak{A}}
\end{equation*}
of the pair $\left( M_{0},Y_{0}\right) $ over $\Delta _{\frak{A}}\subseteq
\Delta $. We realize the deformation $M_{\frak{A}}/\Delta _{\frak{A}}$ as
Kuranishi data $\xi $ on $M_{0}\times \Delta $ associated to a transversely
holomorphic trivialization 
\begin{equation*}
F_{\sigma }:M_{\frak{A}}\rightarrow M_{0}\times \Delta
\end{equation*}
such that 
\begin{equation*}
\overline{\partial }_{M_{0}}\xi -\frac{1}{2}\left[ \xi ,\xi \right] \in
A^{0,2}\left( T_{M_{0}}\right) \otimes \frak{A}
\end{equation*}
and 
\begin{equation*}
Y_{\frak{A}}\subseteq \sigma ^{-1}\left( Y_{0}\right) .
\end{equation*}
The obstruction class for the extension of the pair to a family over $\Delta
_{\frak{mA}}$ is given by the element 
\begin{equation*}
\left( \xi ,\overline{\partial }\xi -\frac{1}{2}\left[ \xi ,\xi \right]
\right) \in \left( A^{0,1}\left( N_{Y_{0}\backslash M_{0}}\right) \oplus
A^{0,2}\left( T_{M_{0}}\right) \right) \otimes \frak{A}
\end{equation*}
in the Dolbeault resolution of the hypercohomology of the complex 
\begin{equation}
T_{M_{0}}\rightarrow N_{Y_{0}\backslash M_{0}}.  \label{5.1.1}
\end{equation}

Indeed, if there is an element 
\begin{equation*}
\left( \alpha ,\beta \right) \in \left( A^{0}\left( N_{Y_{0}\backslash
M_{0}}\right) \oplus A^{0,1}\left( T_{M_{0}}\right) \right) \otimes \frak{A}
\end{equation*}
such that, modulo $\frak{mA}$, 
\begin{equation*}
\delta \left( \alpha ,\beta \right) =\left( -\overline{\partial }\alpha
+\left. \beta \right| _{Y_{0}},\overline{\partial }\beta \right) =\left( \xi
,\ \overline{\partial }\xi -\frac{1}{2}\left[ \xi ,\xi \right] \right) ,
\end{equation*}
then 
\begin{equation*}
\xi -\beta
\end{equation*}
is integrable over $\Delta _{\frak{mA}}$ and, letting $\alpha $ denote a
representative in $A^{0,1}\left( T_{M_{0}}\right) \otimes \frak{A}$ as in \S 
\ref{3}, the trivialization given by the Kuranishi data 
\begin{eqnarray*}
\tilde{\xi} &=&\exp \left( \left[ \alpha ,\ \right] \right) \left( \xi
-\beta \right) -\frac{1-\exp \left( \left[ \alpha ,\ \right] \right) }{%
\left[ \alpha ,\ \right] }\left( \bar{\partial }\alpha \right) \\
&=&\xi -\beta +\left[ \alpha ,\xi -\beta \right] +\bar{\partial }\alpha
+\ldots ,
\end{eqnarray*}
has the property that $\tilde{\sigma}^{-1}\left( Y_{0}\right) $ restricts to
a holomorphic submanifold over $\Delta _{\frak{mA}}$.

Let 
\begin{equation}
T_{Y_{0}|M_{0}}  \label{5.1.2}
\end{equation}
denote the kernel of $\left( \ref{5.1.1}\right) $. Then, since 
\begin{equation*}
\left. \xi \right| _{Y_{0}\times \Delta }\in \left( A^{0,1}\left(
T_{Y_{0}}\right) \otimes \frak{m}+A^{0,1}\left( \left. T_{M_{0}}\right|
_{Y_{0}}\right) \otimes \frak{A}\right) ,
\end{equation*}
we conclude that 
\begin{equation*}
\left. \left[ \xi ,\xi \right] \right| _{Y_{0}\times \Delta }\in
A^{0,1}\left( T_{Y_{0}}\right) \otimes \frak{m}^{2}+A^{0,1}\left( \left.
T_{M_{0}}\right| _{Y_{0}}\right) \otimes \frak{mA}.
\end{equation*}
Let $\varepsilon \in A^{0,1}\left( T_{M_{0}}\right) \otimes \frak{A}$ be
such that 
\begin{equation}
\left( \xi -\varepsilon \right) \in A^{0,1}\left( T_{Y_{0}|M_{0}}\right) .
\label{5.1.20}
\end{equation}
Then 
\begin{equation*}
\left( \overline{\partial }\left( \xi -\varepsilon \right) -\frac{1}{2}
\left[ \xi ,\xi \right] \right) \in \left( A^{0,2}\left(
T_{Y_{0}|M_{0}}\right) \otimes \frak{A}\right)
\end{equation*}
and so $\left\{ \overline{\partial }\left( \xi -\varepsilon \right) -\frac{1%
}{2}\left[ \xi ,\xi \right] \right\} $ in 
\begin{equation*}
H^{2}\left( T_{Y_{0}|M_{0}}\right) \otimes \frac{\frak{A}}{\frak{mA}}\cong 
\Bbb{H}^{2}\left( T_{M_{0}}\rightarrow N_{Y_{0}\backslash M_{0}}\right)
\otimes \frac{\frak{A}}{\frak{mA}}
\end{equation*}
is the obstruction class to extending the deformation of the pair $\left(
M_{0},Y_{0}\right) $ over $\Delta _{\frak{mA}}$.

\begin{theorem}
\label{5.2}Suppose we are given a deformation 
\begin{equation*}
Y_{\frak{A}}\backslash M_{\frak{A}}/\Delta _{\frak{A}}
\end{equation*}
of a compact K\"{a}hler manifold $M_{0}$ and closed submanifold $Y_{0}$ and
obstruction class 
\begin{equation*}
\left\{ \overline{\partial }\left( \xi -\varepsilon \right) -\frac{1}{2}%
\left[ \xi ,\xi \right] \right\} \in H^{2}\left( T_{Y_{0}|M_{0}}\right)
\otimes \frac{\frak{A}}{\frak{mA}}
\end{equation*}
to extending the deformation of the pair over $\Delta _{\frak{mA}}$. For 
\begin{equation*}
\Omega _{Y_{0}|M_{0}}^{p}=\ker \left( \Omega _{M_{0}}^{p}\rightarrow \Omega %
_{Y_{0}}^{p}\right) ,
\end{equation*}
suppose that 
\begin{equation*}
\left\{ \omega _{0}\right\} \in H^{q}\left( \Omega _{Y_{0}|M_{0}}^{p}\right)
=\Bbb{H}^{q}\left( \Omega _{M_{0}}^{p}\rightarrow \Omega _{Y_{0}}^{p}\right)
.
\end{equation*}
Then the value of the pairing 
\begin{equation*}
\begin{array}{ccc}
\left( \left\langle \left. \overline{\partial }\left( \xi -\varepsilon
\right) -\frac{1}{2}\left[ \xi ,\xi \right] \right| \omega _{0}\right\rangle
\right)  & \in  & H^{q+2}\left( \Omega _{Y_{0}|M_{0}}^{p-1}\right) \otimes 
\frac{\frak{A}}{\frak{mA}} \\ 
&  & = \\ 
&  & \Bbb{H}^{q+2}\left( \Omega _{M_{0}}^{p-1}\rightarrow \Omega %
_{Y_{0}}^{p-1}\right) \otimes \frac{\frak{A}}{\frak{mA}}
\end{array}
\end{equation*}
is zero.
\end{theorem}

\begin{proof}
By assumption we have at the level of cohomology that 
\begin{equation*}
\left. \omega _{0}\right| _{Y_{0}}\equiv 0.
\end{equation*}
As in the proof of Theorem \ref{4.1}, there is an extension of $\omega _{0}$
to a $\partial $-closed form $\omega $ of type $\left( p,q\right) $ on $%
M_{0}\times \Delta $ such that 
\begin{equation*}
\left. \left( \overline{\partial }-L^{1,0}\right) \omega \right| _{M_{\frak{A%
}}}=0.
\end{equation*}
By hypothesis 
\begin{equation*}
\left. \varphi ^{-1}\left( e^{\left\langle \left. \xi \right| \
\right\rangle }\omega \right) \right| _{Y_{\frak{A}}}
\end{equation*}
is $\overline{\partial }_{Y_{\frak{A}}}$-exact, that is, at the level of
cohomology we have 
\begin{equation}
\left. \varphi ^{-1}\left( e^{\left\langle \left. \xi \right| \
\right\rangle }\omega \right) \right| _{Y_{\frak{A}}}\equiv 0.  \label{5.2.2}
\end{equation}
So by the Hodge $\partial \overline{\partial }$-lemma on $M_{\frak{A}}$
there is 
\begin{equation*}
\gamma \in B^{p-1,q-1}\left( M/\Delta \right)
\end{equation*}
such that 
\begin{equation*}
\left. \partial _{M}\overline{\partial }_{M}\gamma \right| _{Y_{\frak{A}
}}=\left. \varphi ^{-1}e^{\left\langle \left. \xi \right| \ \right\rangle
}\omega \right| _{Y_{\frak{A}}}.
\end{equation*}
Adjusting the choice of $\omega $ by $\partial _{M}\overline{\partial }
_{M}\gamma $ we can assume that 
\begin{equation*}
\left. \varphi ^{-1}e^{\left\langle \left. \xi \right| \ \right\rangle
}\omega \right| _{Y_{\frak{A}}}=0
\end{equation*}
as a cocycle, that is 
\begin{equation*}
\left. \omega \right| _{Y_{0}\times \Delta _{\frak{A}}}=0.
\end{equation*}

Referring to $\left( \ref{5.1.20}\right) $ let 
\begin{equation*}
\xi ^{\prime }=\xi -\varepsilon 
\end{equation*}
and notice that 
\begin{equation}
\left[ \xi ^{\prime },\xi ^{\prime }\right] \equiv \left[ \xi ,\xi \right] 
\label{5.2.21}
\end{equation}
modulo $\frak{mA}$. As before by $\left( \ref{2.7.7}\right) $ and the fact
that $\omega $ is $\partial $-closed, we have 
\begin{eqnarray*}
-\left\langle \left. \left[ \xi ^{\prime },\xi ^{\prime }\right] \right|
\omega \right\rangle  &=&\left\langle \left. \xi ^{\prime }\right| L_{\xi
^{\prime }}^{1,0}\omega \right\rangle -L_{\xi ^{\prime }}^{1,0}\left\langle
\left. \xi ^{\prime }\right| \omega \right\rangle  \\
&=&-2\left\langle \left. \xi ^{\prime }\right| \partial \left\langle \left.
\xi ^{\prime }\right| \omega \right\rangle \right\rangle +\partial %
\left\langle \left. \xi ^{\prime }\right| \left\langle \left. \xi ^{\prime
}\right| \omega \right\rangle \right\rangle 
\end{eqnarray*}
so that 
\begin{eqnarray}
\frac{1}{2}\left\langle \left. \left[ \xi ^{\prime },\xi ^{\prime }\right]
\right| \omega \right\rangle  &=&\left\langle \left. \xi ^{\prime }\right| 
\partial \left\langle \left. \xi ^{\prime }\right| \omega \right\rangle
\right\rangle -\frac{1}{2}\partial \left\langle \left. \xi ^{\prime }\right|
\left\langle \left. \xi ^{\prime }\right| \omega \right\rangle \right\rangle 
\label{5.2.3} \\
&=&-\left\langle \left. \xi ^{\prime }\right| L_{\xi ^{\prime }}\omega
\right\rangle -\frac{1}{2}\partial \left\langle \left. \xi ^{\prime }\right|
\left\langle \left. \xi ^{\prime }\right| \omega \right\rangle \right\rangle
.  \notag
\end{eqnarray}
Since 
\begin{equation*}
\left( \overline{\partial }_{M_{0}}-L_{\xi ^{\prime }}^{1,0}\right) \omega
\in A^{*}\left( M_{0}\right) \otimes \frak{A},
\end{equation*}
and 
\begin{equation*}
\xi ^{\prime }\in A^{0,1}\left( T_{M_{0}}\right) \otimes \frak{m},
\end{equation*}
we can compute modulo $\frak{mA}$ as in the proof of Theorem \ref{4.1}: 
\begin{eqnarray*}
\overline{\partial }\left\langle \left. \xi ^{\prime }\right| \omega
\right\rangle  &=&\left\langle \left. \overline{\partial }\xi ^{\prime
}\right| \omega \right\rangle +\left\langle \left. \xi ^{\prime }\right| 
\overline{\partial }\omega \right\rangle  \\
&\equiv &\left\langle \left. \overline{\partial }\xi ^{\prime }\right|
\omega \right\rangle +\left\langle \left. \xi ^{\prime }\right| L_{\xi
^{\prime }}^{1,0}\omega \right\rangle  \\
&=&\left\langle \left. \overline{\partial }\xi ^{\prime }\right| \omega
\right\rangle -\frac{1}{2}\left\langle \left. \left[ \xi ^{\prime },\xi
^{\prime }\right] \right| \omega \right\rangle -\frac{1}{2}\partial %
\left\langle \left. \xi ^{\prime }\right| \left\langle \left. \xi ^{\prime
}\right| \omega \right\rangle \right\rangle  \\
&\equiv &\left\langle \left. \overline{\partial }\xi ^{\prime }-\frac{1}{2}%
\left[ \xi ^{\prime },\xi ^{\prime }\right] \right| \omega _{0}\right\rangle
-\frac{1}{2}\partial \left\langle \left. \xi ^{\prime }\right| \left\langle
\left. \xi ^{\prime }\right| \omega \right\rangle \right\rangle .
\end{eqnarray*}
So by $\left( \ref{5.2.21}\right) $ we have modulo $\frak{mA}$ that 
\begin{equation*}
\overline{\partial }\left\langle \left. \xi ^{\prime }\right| \omega
\right\rangle +\frac{1}{2}\partial \left\langle \left. \xi ^{\prime }\right|
\left\langle \left. \xi ^{\prime }\right| \omega _{0}\right\rangle
\right\rangle \equiv \overline{\partial }\left\langle \left. \xi ^{\prime
}\right| \omega _{0}\right\rangle -\frac{1}{2}\left\langle \left. \left[ \xi
,\xi \right] \right| \omega _{0}\right\rangle 
\end{equation*}
and so 
\begin{equation*}
\left\langle \left. \left[ \xi ,\xi \right] \right| \omega _{0}\right\rangle
\equiv 0\in H^{p-1,q+2}\left( M_{0},Y_{0}\right) \otimes \frac{\Bbb{C}\left[
\left[ t\right] \right] }{\frak{mA}}.
\end{equation*}
\end{proof}

Notice that, while Theorem \ref{5.2} fully generalizes Theorem \ref{4.1}, it
does not quite generalize the Semiregularity Theorem \ref{4.3} since, under
the hypotheses of Theorem \ref{4.3}, Theorem \ref{5.2} only yields the
weaker result that 
\begin{equation*}
\left\langle \left. \xi \right| \omega _{0}\right\rangle \in image\left(
H^{p-1,q+1}\left( M_{0}\right) \rightarrow H^{p-1,q+1}\left( Y_{0}\right)
\right) \otimes \frac{\frak{A}}{\frak{mA}}.
\end{equation*}
We can however fully generalize the Semiregularity Theorem as well. To do
this we let 
\begin{equation*}
A_{Y_{0}}^{0,k}\left( T_{M_{0}}\right)
\end{equation*}
be the subspace of $A^{0,k}\left( T_{M_{0}}\right) $ consisting of those $%
\xi _{0}$ such that 
\begin{equation*}
\left. \left\langle \left. \xi _{0}\right| \eta _{0}\right\rangle \right|
_{Y_{0}}
\end{equation*}
is $d$-exact for every $d$-closed form $\eta _{0}$ on $M_{0}$ such that $%
\left. \eta _{0}\right| _{Y_{0}}$ is $d$-exact. By $\left( \ref{2.7.8}
\right) $ $A_{Y_{0}}^{0,*}\left( T_{M_{0}}\right) \otimes \Bbb{C}\left[
\left[ t\right] \right] $ is a differential graded Lie sub-algebra of $%
A^{0,*}\left( T_{M_{0}}\right) \otimes \Bbb{C}\left[ \left[ t\right] \right]
.$ Then $A_{Y_{0}}^{0,*}\left( T_{M_{0}}\right) \otimes \Bbb{C}\left[ \left[
t\right] \right] $ is the differential graded Lie algebra which measures the
deformations of $M_{0}$ such that 
\begin{equation*}
K:=\ker \left( H^{*}\left( M_{0}\right) \rightarrow H^{*}\left( Y_{0}\right)
\right)
\end{equation*}
is a rational sub-Hodge structure$.$ It is then a tautology that
obstructions in 
\begin{equation*}
H^{1}\left( A_{Y_{0}}^{0,*}\left( N_{Y_{0}|M_{0}}\right) ,\overline{\partial 
}\right)
\end{equation*}
annihilate $K$.

\part{Appendices\label{part5}}

\section{Existence of transversely holomorphic trivializations\label{A1}}

\begin{lemma}
\label{2.3}i) There exist transversely holomorphic trivializations of any
deformation $\left( \ref{2.1}\right) $.

ii) Each transversely holomorphic trivialization 
\begin{equation*}
F:M\overset{\left( \sigma ,\ \right) }{\longrightarrow }M_{0}\times \Delta
\end{equation*}
of a deformation $M/\Delta $ as in Definition \ref{2.2} determines Kuranshi
data 
\begin{equation*}
\xi _{\sigma }:\Delta \rightarrow A_{M_{0}}^{0,1}\left( T_{1,0}\right) .
\end{equation*}
iii) If two deformations/transversely holomorphic trivializations 
\begin{eqnarray*}
F &:&M\overset{\left( \sigma ,\ \right) }{\longrightarrow }M_{0}\times \Delta
\\
F^{\prime } &:&M^{\prime }\overset{\left( \sigma ^{\prime },\ \right) }{%
\longrightarrow }M_{0}\times \Delta
\end{eqnarray*}
induce the same Kuranishi data 
\begin{equation*}
\xi _{\sigma }=\xi _{\sigma ^{\prime }}:\Delta \rightarrow
A_{M_{0}}^{0,1}\left( T_{1,0}\right) ,
\end{equation*}
then there is a holomorphic isomorphism 
\begin{equation*}
\varphi :M\rightarrow M^{\prime }
\end{equation*}
defined over $\Delta $ such that the diagram 
\begin{equation*}
\begin{array}{lllll}
M &  & \overset{\varphi }{\longrightarrow } &  & M^{\prime } \\ 
& \searrow ^{F} &  & \swarrow _{F^{\prime }} &  \\ 
&  & M_{0}\times \Delta &  & 
\end{array}
\end{equation*}
commutes.

iv) Each transversely holomorphic trivialization as in Definition \ref{2.2}
induces a unique lifting 
\begin{equation*}
\tau _{j}
\end{equation*}
of $\frac{\partial }{\partial t_{j}}$ to a $C^{\infty }$-vector field of
type $\left( 1,0\right) $ on $M$ and all liftings occur in some
trivialization of $M/\Delta $.
\end{lemma}

\begin{proof}
Once we construct a transversely holomorphic trivialization $\sigma $ of $%
M/\Delta $, it will suffice to establish that, given a point $x_{0}\in M_{0}$
and a local holomorphic coordinate system 
\begin{equation*}
v_{W_{0}}
\end{equation*}
on a neighborhood $W_{0}$ of $x_{0}$ in $M_{0}$, there exists a local
holomorphic coordinate system 
\begin{equation*}
\left( v_{W},t\right) =\left( v_{W,x_{0}},t\right)
\end{equation*}
on a (formal) neighborhood $W$ of $W_{0}$ with the four properties below:

\begin{enumerate}
\item  \label{2.3.1}At points of $\sigma ^{-1}\left( x_{0}\right) $, 
\begin{equation}
dv_{W}=\left. \sigma ^{*}dv_{W_{0}}\left( x\right) \right|
_{x=x_{0}}+\sum\nolimits_{J}\left. t^{J}C_{J}\sigma ^{*}d\overline{v_{W_{0}}}%
\left( x\right) \right| _{x=x_{0}}  \label{formexp}
\end{equation}
for some system of $m\times m$ matrices $C_{J}=\left( c_{J,s}^{r}\right) $.
(Here and in what follows all differentials are relative over $\Delta $.)
Thus, modulo the ideal 
\begin{equation*}
\overline{t}=\left\{ \overline{t}_{j}\right\} ,
\end{equation*}
we also have 
\begin{eqnarray*}
d\overline{v_{W}^{s}} &\equiv &d\overline{v_{W_{0}}^{s}} \\
\frac{\partial }{\partial \overline{v_{W}^{s}}} &=&\frac{\partial }{\partial 
\overline{v_{W_{0}}^{s}}}-\sum\nolimits_{J,\left| J\right|
>0}c_{I,s}^{r}\left( x_{0}\right) t^{J}\frac{\partial }{\partial
v_{W_{0}}^{r}} \\
\frac{\partial }{\partial v_{W}^{s}} &\equiv &\frac{\partial }{\partial
v_{W_{0}}^{s}}.
\end{eqnarray*}

\item  \label{2.3.2}(Using the Einstein summation convention) the mapping 
\begin{equation}
\xi _{\sigma }=\sum\nolimits_{\left| J\right| >0}\xi _{J}t^{J}
\label{convser}
\end{equation}
is given by the tensor 
\begin{equation*}
\xi _{J}:=d\overline{v_{W_{0}}^{k}}\otimes c_{J,k}^{l}\frac{\partial }{%
\partial v_{W_{0}}^{l}}\in A^{0,1}\left( T_{M_{0}}^{1,0}\right)
\end{equation*}
which is independent of the choice of coordinates satisfying $\ref{2.3.1}$
and $\ref{2.3.2}$. That is, if $\hat{v}_{W}$ is another choice of
holomorphic coordinates satisfying $\ref{2.3.1}$ and $\ref{2.3.2}$ along $%
\sigma ^{-1}\left( x_{0}\right) $, then at points of $\sigma ^{-1}\left(
x_{0}\right) ,$%
\begin{equation*}
d\overline{v_{W_{0}}^{k}}\otimes \left\langle \left. \sum\nolimits_{\left|
J\right| >0}c_{J,k}^{l}t^{J}\frac{\partial }{\partial v_{W_{0}}^{l}}\right|
\;\right\rangle =d\overline{\hat{v}_{W_{0}}^{k}}\otimes \left\langle \left.
\sum\nolimits_{\left| J\right| >0}\hat{c}_{J,k}^{l}t^{J}\frac{\partial }{%
\partial \hat{v}_{W_{0}}^{l}}\right| \;\right\rangle .
\end{equation*}
Notice that the series $\left( \ref{convser}\right) $ is convergent since
locally it equals $\xi _{\beta }$ for a convergent series 
\begin{equation*}
\beta =\sum\nolimits_{\left| J\right| >0}\beta _{J}t^{J}.
\end{equation*}
Thus the formal expression $\left( \ref{formexp}\right) $ is also convergent.

\item  \label{2.3.3}If 
\begin{equation*}
\xi _{\sigma }=\xi _{\sigma ^{\prime }}
\end{equation*}
for 
\begin{eqnarray*}
F &:&M\overset{\left( \sigma ,\ \right) }{\longrightarrow }M_{0}\times \Delta
\\
F^{\prime } &:&M^{\prime }\overset{\left( \sigma ^{\prime },\ \right) }{%
\longrightarrow }M_{0}\times \Delta ,
\end{eqnarray*}
then 
\begin{equation*}
\left( F^{\prime }\right) ^{-1}\circ F
\end{equation*}
is holomorphic.

\item  \label{2.3.4}Each transversely holomorphic trivialization as in
Definition\ref{2.2} induces a unique lifting 
\begin{equation*}
\tau _{j}
\end{equation*}
of $\frac{\partial }{\partial t_{j}}$ to a $C^{\infty }$-vector field of
type $\left( 1,0\right) $ on $M$. Given any $C^{\infty }$-section $\vartheta
\in A^{0,0}\left( T_{1,0}\left( M_{0}\right) \right) $, there is another
transversely holomorphic trivialization $\sigma ^{\prime }$ such that 
\begin{equation*}
\left. \tau _{j}^{\prime }\right| _{M_{0}}=\left. \tau _{j}\right|
_{M_{0}}+\vartheta .
\end{equation*}
\end{enumerate}

To begin we must construct a transversely holomorphic trivialization of $%
M/\Delta $. Consider the Grassmann bundle of linear subspaces of $\left.
T_{1,0}\left( M\right) \right| _{M_{0}}$ of dimension complementary to that
of $T_{1,0}\left( M_{0}\right) .$ At each point $x_{0}$ of $M_{0}$, the set
of subspaces transverse to $\left. T_{1,0}\left( M_{0}\right) \right|
_{x_{0}}$ form an affine space. Thus the set of sections $V$ of 
\begin{equation*}
Gr\left( \dim M_{0,}\left. T_{1,0}\left( M\right) \right| _{M_{0}}\right)
\end{equation*}
which, together with $\Bbb{P}\left( T_{1,0}\left( M_{0}\right) \right) $,
generates the fiber of the projective bundle $\Bbb{P}\left( \left.
T_{1,0}\left( M\right) \right| _{M_{0}}\right) $ at each point of $M_{0}$ is
convex. Thus the set of choices is path connected. For each choice, the
morphism 
\begin{equation*}
V\rightarrow T_{1,0}\left( \Delta \right) _{0}
\end{equation*}
gives a distinguished framing. One can choose a covering of $M_{0}$ by small
open sets $W$ in $M$, each of which have a coordinate system $\left(
u_{W},t\right) $, such that the sets 
\begin{equation*}
u_{W}=constant
\end{equation*}
have tangent space which approximates the subspace of $\Bbb{P}\left( \left.
T_{1,0}\left( M\right) \right| _{M_{0}}\right) $ spanned by the chosen
sections. Next choose a $C^{\infty }$-partition-of-unity $\left\{ \rho
_{W}\right\} $ on $M_{0}$ subordinate to the covering $\left\{
W_{0}:=M_{0}\cap W\right\} $. Let $\Gamma $ denote the diagonal of $%
M_{0}\times M_{0}$ considered as a submanifold of $M\times M_{0}.$ The graph
of the projection in Definition \ref{2.2} is then given, for $\left(
y,x\right) \in W\times W_{0}$ (and $y$ sufficiently near $x$), by the
equation 
\begin{equation}
\sum\nolimits_{W^{\prime }}\rho _{W^{\prime }}\left( x\right) \frac{\partial
u_{W}}{\partial u_{W^{\prime }}}\left( x\right) \left( u_{W^{\prime }}\left(
y\right) -u_{W^{\prime }}\left( x\right) \right) =0.  \label{2.4.1}
\end{equation}
To see that the trivialization is well-defined, independently of the choice
of the coordinate patch $W$, notice that, for a second coordinate patch $V$,
we have 
\begin{eqnarray*}
&&\frac{\partial u_{V}}{\partial u_{W}}\left( x\right) \left(
\sum\nolimits_{W^{\prime }}\rho _{W^{\prime }}\left( x\right) \frac{\partial
u_{W}}{\partial u_{W^{\prime }}}\left( x\right) \left( u_{W^{\prime }}\left(
y\right) -u_{W^{\prime }}\left( x\right) \right) \right) \\
&=&\left( \sum\nolimits_{W^{\prime }}\rho _{W^{\prime }}\left( x\right) 
\frac{\partial u_{V}}{\partial u_{W}}\left( x\right) \cdot \frac{\partial
u_{W}}{\partial u_{W^{\prime }}}\left( x\right) \left( u_{W^{\prime }}\left(
y\right) -u_{W^{\prime }}\left( x\right) \right) \right) \\
&=&\sum\nolimits_{W^{\prime }}\rho _{W^{\prime }}\left( x\right) \frac{
\partial u_{V}}{\partial u_{W^{\prime }}}\left( x\right) \left( u_{W^{\prime
}}\left( y\right) -u_{W^{\prime }}\left( x\right) \right) .
\end{eqnarray*}
Next fix a point $x_{0}\in M_{0}$. Use holomorphic local coordinate 
\begin{equation}
\tilde{v}_{W}\left( y\right) =\sum\nolimits_{W^{\prime }}\rho _{W^{\prime
}}\left( x_{0}\right) \frac{\partial u_{W}}{\partial u_{W^{\prime }}}\left(
x_{0}\right) \left( u_{W^{\prime }}\left( y\right) -u_{W^{\prime }}\left(
x_{0}\right) \right)  \label{2.4.2}
\end{equation}
on $W$. At each point $y\in \sigma ^{-1}\left( x_{0}\right) $ we have, by $%
\left( \ref{2.4.1}\right) $, that 
\begin{equation*}
\begin{tabular}{l}
$d\tilde{v}_{W}\left( y\right) =\sum\nolimits_{W^{\prime }}\rho _{W^{\prime
}}\left( x_{0}\right) \frac{\partial u_{W}}{\partial u_{W^{\prime }}}\left(
x_{0}\right) du_{W^{\prime }}\left( y\right) $ \\ 
$=\sum\nolimits_{W^{\prime }}\rho _{W^{\prime }}\left( x_{0}\right) \frac{
\partial u_{W}}{\partial u_{W^{\prime }}}\left( x_{0}\right) \left.
du_{W^{\prime }}\left( x\right) \right| _{x=x_{0}}$ \\ 
\quad \quad $-\sum\nolimits_{W^{\prime }}\left( u_{W^{\prime }}\left(
y\right) -u_{W^{\prime }}\left( x_{0}\right) \right) \left. d\left( \rho
_{W^{\prime }}\left( x\right) \frac{\partial u_{W}}{\partial u_{W^{\prime }}}
\left( x\right) \right) \right| _{x=x_{0}}$ \\ 
$=\sum\nolimits_{W^{\prime }}\rho _{W^{\prime }}\left( x_{0}\right) \frac{
\partial u_{W}}{\partial u_{W^{\prime }}}\left( x_{0}\right) \left.
du_{W^{\prime }}\left( x\right) \right| _{x=x_{0}}$ \\ 
\quad \quad $-\sum\nolimits_{W^{\prime }}\sum\nolimits_{I}a_{I}t^{I}\left.
d\left( \rho _{W^{\prime }}\left( x\right) \frac{\partial u_{W}}{\partial
u_{W^{\prime }}}\left( x\right) \right) \right| _{x=x_{0}}$ \\ 
$=\left. d\tilde{v}_{W_{0}}\left( x\right) \right| _{x=x_{0}}$ \\ 
\quad \quad $+\sum\nolimits_{J,\left| J\right| >0}t^{J}\left. A_{J}d\tilde{v}
_{W_{0}}\left( x\right) \right| _{x=x_{0}}$ \\ 
\quad \quad \quad \quad $+\sum\nolimits_{J,\left| J\right| >0}t^{J}\left.
B_{J}d\overline{\tilde{v}_{W_{0}}}\left( x\right) \right| _{x=x_{0}}$%
\end{tabular}
\end{equation*}
for some systems of $m\times m$ matrices $A_{J}$ and $B_{J}$.

Suppose inductively that, for $\left| J\right| <n$, 
\begin{equation*}
A_{J}=0.
\end{equation*}
For new holomorphic coordinates 
\begin{equation*}
\tilde{v}_{W}^{\prime }=\tilde{v}_{W}-\sum t^{J}A_{J}\tilde{v}_{W}
\end{equation*}
on $W$ we have: 
\begin{eqnarray*}
d\tilde{v}_{W}^{\prime } &=& \\
&&\left. d\tilde{v}_{W_{0}}\left( x\right) \right| _{x=x_{0}} \\
&&\quad \quad +\sum\nolimits_{\left| J\right| >n}t^{J}\left. A_{J}d\tilde{v}
_{W_{0}}\left( x\right) \right| _{x=x_{0}} \\
&&\quad \quad \quad \quad +\sum\nolimits_{J,\left| J\right| >0}t^{J}\left.
B_{J}^{\prime }d\overline{\tilde{v}_{W_{0}}}\left( x\right) \right|
_{x=x_{0}}.
\end{eqnarray*}
Repeating this construction, we have at least a formal set of holomorphic
coordinates $v_{W}$ on $W$ such that at each point $y\in \sigma ^{-1}\left(
x_{0}\right) $%
\begin{equation}
dv_{W}=dv_{W_{0}}\left( x_{0}\right) +\sum\nolimits_{J,\left| J\right|
>0}^{\infty }t^{J}C_{J}d\overline{v_{W_{0}}}\left( x_{0}\right) .
\label{2.4.2.5}
\end{equation}
where this series in $t$ with values in $\left. T_{1,0}\left( M_{0}\right)
\right| _{x_{0}}\oplus \left. T_{0,1}\left( M_{0}\right) \right| _{x_{0}}$
is convergent. That is, in the cotangent space at points $y\in \sigma
^{-1}\left( x_{0}\right) $, the subspace of one-forms of type $\left(
1,0\right) $ is exactly the subspace annihilated by the vectors 
\begin{equation}
\frac{\partial }{\partial \overline{v_{W_{0}}^{s}}}-\sum\nolimits_{J,\left|
J\right| >0}c_{J,s}^{r}\left( x_{0}\right) t^{J}\frac{\partial }{\partial
v_{W_{0}}^{r}}  \label{2.4.3}
\end{equation}
where each $C_{J}=\left( c_{J,s}^{r}\right) $ is some $C^{\infty }$
-matrix-valued function of $x_{0}$. This proves \ref{2.3.1}.

To prove \ref{2.3.2} let 
\begin{equation*}
\varphi :W/\Delta \rightarrow W/\Delta
\end{equation*}
be the holomorphic automorphism such that 
\begin{equation*}
\hat{v}_{W}=v_{W}\circ \varphi
\end{equation*}
and define the trivialization $\hat{\sigma}$ by commutativity of the diagram 
\begin{equation*}
\begin{array}{ccc}
W & \overset{\varphi }{\rightarrow } & W \\ 
\downarrow \hat{\sigma} &  & \downarrow \sigma \\ 
W_{0} & \overset{\varphi _{0}}{\rightarrow } & W_{0}
\end{array}
\end{equation*}
where $\varphi _{0}=\left. \varphi \right| _{W_{0}}$. Again computing modulo 
$\overline{t}$ we write 
\begin{equation*}
\overline{\partial }\equiv d\overline{v_{W}^{s}}\otimes \frac{\partial }{
\partial \overline{v_{W}^{s}}}\equiv d\overline{v_{W_{0}}^{s}}\otimes \left( 
\frac{\partial }{\partial \overline{v_{W_{0}}^{s}}}-\sum\nolimits_{J,\left|
J\right| >0}c_{J,s}^{r}\left( x_{0}\right) t^{J}\frac{\partial }{\partial
v_{W_{0}}^{r}}\right) .
\end{equation*}
But since $\varphi $ is holomorphic 
\begin{equation*}
\overline{\partial }=\varphi ^{*}\overline{\partial }\equiv d\overline{\hat{v%
}_{W_{0}}^{s}}\otimes \left( \frac{\partial }{\partial \overline{\hat{v}
_{W_{0}}^{s}}}-\sum\nolimits_{J,\left| J\right| >0}\hat{c}_{J,s}^{r}\left(
x_{0}\right) t^{J}\frac{\partial }{\partial \hat{v}_{W_{0}}^{r}}\right) .
\end{equation*}

To prove \ref{2.3.3} notice that, at each point of $\sigma ^{-1}\left(
x_{0}\right) $ the space 
\begin{equation*}
\left. T^{1,0}\left( M/\Delta \right) \right| _{\sigma ^{-1}\left(
x_{0}\right) }
\end{equation*}
corresponds under $F_{*}$ to the family of subspaces annihilated by 
\begin{equation*}
\frac{\partial }{\partial \overline{v_{W_{0}}^{s}}}-\sum\nolimits_{J,\left|
J\right| >0}c_{J,s}^{r}\left( x_{0}\right) t^{J}\frac{\partial }{\partial
v_{W_{0}}^{r}}
\end{equation*}
for all $s$ and similarly for $F^{\prime }$. Thus for all $x_{0}\in M_{0}$
we have 
\begin{equation*}
\left( \left( F^{\prime }\right) ^{-1}\circ F\right) ^{*}\left( \left.
T^{1,0}\left( M/\Delta \right) \right| _{\left( \sigma ^{\prime }\right)
^{-1}\left( x_{0}\right) }\right) \subseteq \left. T^{1,0}\left( M/\Delta
\right) \right| _{\sigma ^{-1}\left( x_{0}\right) }
\end{equation*}
so that $\left( F^{\prime }\right) ^{-1}\circ F$ is holomorphic.

To prove \ref{2.3.4} write 
\begin{equation*}
\tilde{v}_{W}^{\prime }=v_{W}+\vartheta \left( v_{W_{0}}\right) t_{j}
\end{equation*}
and repeat the normalization process as above from $\left( \ref{2.4.2}
\right) $ forward with $\tilde{v}_{W}^{\prime }$ replacing $\tilde{v}_{W}$.
\end{proof}

\section{Lie derivatives, standard identities}

We make precise the two actions of an element $\xi \in A^{0,k}\left(
T_{1,0}\left( M_{0}\right) \right) $ on $\sum A^{p,q}\left( M_{0}\right) $
and review the elementary identities for these actions that are used in this
paper. We write the action via contraction as 
\begin{equation}
\left\langle \xi \right. \left| \ \right\rangle ,  \label{2.7.0}
\end{equation}
and ``Lie differentiation'' as 
\begin{equation}
L_{\xi }:=\left\langle \xi \right. \left| \ \right\rangle \circ d+\left(
-1\right) ^{k}d\circ \left\langle \xi \right. \left| \ \right\rangle .
\label{2.7.1}
\end{equation}
(See also \S 5.3 of \cite{Ko}.) The sign is so chosen that, writing any
element of $A^{0,k}\left( T_{1,0}\left( M_{0}\right) \right) $ locally as a
sum of terms 
\begin{equation*}
\xi =\bar{\eta}\otimes \chi
\end{equation*}
for some $d$-closed $\left( 0,k\right) $ -form $\bar{\eta}$ and $\chi \in
A^{0,0}\left( T_{1,0}\left( M_{0}\right) \right) $, then 
\begin{equation*}
L_{\xi }=\bar{\eta}\otimes L_{\chi }.
\end{equation*}
Warning: Since 
\begin{equation*}
\overline{\partial }\left\langle \xi \right. \left| \ \right\rangle
=\left\langle \overline{\partial }\xi \right. \left| \ \right\rangle +\left(
-1\right) ^{k+1}\left\langle \xi \right. \left| \overline{\partial }\
\right\rangle
\end{equation*}
we only obtain that for $k$ \textit{odd} do we have 
\begin{equation*}
\left[ \overline{\partial },L_{\xi }\right] =L_{\overline{\partial }\xi
}:A^{p,q}\left( M_{0}\right) \rightarrow A^{p,q+k+1}\left( M_{0}\right)
+A^{p-1,q+k+2}.
\end{equation*}
However $\left[ \overline{\partial },L_{\xi }\right] $ and $L_{\overline{%
\partial }\xi }$ always act as the same operator on $A^{0,q}\left(
M_{0}\right) $.

We decompose 
\begin{equation*}
L_{\xi }:A^{p,q}\left( M_{0}\right) \rightarrow A^{p,q+k}\left( M_{0}\right)
+A^{p-1,q+k+1}\left( M_{0}\right)
\end{equation*}
into 
\begin{eqnarray}
L_{\xi }^{1,0} &:&=\left\langle \xi \right. \left| \ \right\rangle \circ
\partial +\left( -1\right) ^{k}\partial \circ \left\langle \xi \right.
\left| \ \right\rangle  \label{2.7.10} \\
L_{\xi }^{0,1} &:&=\left\langle \xi \right. \left| \ \right\rangle \circ 
\overline{\partial }+\left( -1\right) ^{k}\overline{\partial }\circ
\left\langle \xi \right. \left| \ \right\rangle .  \notag
\end{eqnarray}

We will review the commutater properties of the operators $\left( \ref{2.7.1}
\right) $ below. For the $\mathcal{O}_{M_{0}}$-linear operators $\left( \ref
{2.7.0}\right) ,$ the situation is much easier, namely, 
\begin{equation*}
\left\langle \xi _{1}\right. \left| \ \right\rangle \circ \left\langle \xi
_{2}\right. \left| \ \right\rangle +\left( -1\right) ^{k}\left\langle \xi
_{2}\right. \left| \ \right\rangle \circ \left\langle \xi _{1}\right. \left|
\ \right\rangle =0
\end{equation*}
so that, in particular, for $k=1$, 
\begin{equation}
\left[ \left\langle \xi _{1}\right. \left| \ \right\rangle ,\left\langle \xi
_{2}\right. \left| \ \right\rangle \right] =0.  \label{2.7.00}
\end{equation}

Also we compute 
\begin{eqnarray}
L_{\xi }L_{\xi ^{\prime }}-\left( -1\right) ^{\deg \bar{\eta}\cdot \deg \bar{
\eta}^{\prime }}L_{\xi ^{\prime }}L_{\xi } &=&\left( \bar{\eta}\otimes
L_{\chi }\right) \left( \bar{\eta}^{\prime }\otimes L_{\chi ^{\prime
}}\right)  \label{48.5} \\
&&-\left( -1\right) ^{\deg \bar{\eta}\cdot \deg \bar{\eta}^{\prime }}\left( 
\bar{\eta}^{\prime }\otimes L_{\chi ^{\prime }}\right) \left( \bar{\eta}
\otimes L_{\chi }\right)  \notag \\
&=&\bar{\eta}\bar{\eta}^{\prime }\left( L_{\chi }L_{\chi ^{\prime }}-L_{\chi
^{\prime }}L_{\chi }\right)  \notag \\
&=&\bar{\eta}\bar{\eta}^{\prime }L_{\left[ \chi ,\chi ^{\prime }\right] }. 
\notag
\end{eqnarray}
So, using this local presentation for 
\begin{equation*}
\xi \in A^{0,j}\left( T_{1,0}\left( M_{0}\right) \right) ,\xi ^{\prime }\in
A^{0,k}\left( T_{1,0}\left( M_{0}\right) \right) ,
\end{equation*}
we can define 
\begin{equation}
\left[ \xi ,\xi ^{\prime }\right] =\bar{\eta}\bar{\eta}^{\prime }\left[ \chi
,\chi ^{\prime }\right] \in A^{0,j+k}\left( T_{1,0}\left( M_{0}\right)
\right) .  \label{2.7.2}
\end{equation}
Notice also that we have the formula 
\begin{eqnarray*}
L_{\xi }\left( \omega \wedge \eta \right) &=&\left\langle \xi \right. \left|
d\omega \wedge \eta +\left( -1\right) ^{\deg \omega }\omega \wedge d\eta
\right\rangle \\
&&+\left( -1\right) ^{k}d\left( \left\langle \xi \right. \left| \omega
\right\rangle \wedge \eta +\left( -1\right) ^{\left( k+1\right) \deg \omega
}\omega \wedge \left\langle \xi \right. \left| \eta \right\rangle \right) \\
&=&\left\langle \xi \right. \left| d\omega \right\rangle \wedge \eta +\left(
-1\right) ^{\left( k+1\right) \left( \deg \omega +1\right) }d\omega \wedge
\left\langle \xi \right. \left| \eta \right\rangle \\
&&+\left( -1\right) ^{\deg \omega }\left\langle \xi \right. \left| \omega
\right\rangle \wedge d\eta +\left( -1\right) ^{k\deg \omega }\omega \wedge
\left\langle \xi \right. \left| d\eta \right\rangle \\
&&+\left( -1\right) ^{k}\left( 
\begin{array}{c}
\left( d\left\langle \xi \right. \left| \omega \right\rangle \right) \wedge
\eta +\left( -1\right) ^{\deg \omega -1+k}\left( \left\langle \xi \right.
\left| \omega \right\rangle \wedge d\eta \right) \\ 
+\left( -1\right) ^{\left( k+1\right) \deg \omega }\left( d\omega \wedge
\left\langle \xi \right. \left| \eta \right\rangle \right) +\left( -1\right)
^{k\deg \omega }\left( \omega \wedge d\left\langle \xi \right. \left| \eta
\right\rangle \right)
\end{array}
\right) \\
&=&L_{\xi }\left( \omega \right) \wedge \eta +\left( -1\right) ^{k\deg
\omega }\omega \wedge L_{\xi }\left( \eta \right) \\
&&+\left( \left( -1\right) ^{\deg \omega }+\left( -1\right) ^{k\left( \deg
\omega -1\right) }\right) \left( \left\langle \xi \right. \left| \omega
\right\rangle \wedge d\eta \right)
\end{eqnarray*}
so that, for $\deg \xi =1$, 
\begin{equation}
L_{\xi }\left( \omega \wedge \eta \right) =L_{\xi }\left( \omega \right)
\wedge \eta +\left( -1\right) ^{\deg \omega }\omega \wedge L_{\xi }\left(
\eta \right) .  \label{2.7.3}
\end{equation}
We have the follow additional general formulas for arbitrary degree: 
\begin{eqnarray}
\bar{\partial }\left[ \xi ,\xi ^{\prime }\right] &=&\bar{\partial }\bar{\eta}%
\bar{\eta}^{\prime }\left[ \chi ,\chi ^{\prime }\right] =\left( -1\right)
^{\deg \bar{\eta}\bar{\eta}^{\prime }}\bar{\eta}\bar{\eta}^{\prime }\bar{
\partial }\left[ \chi ,\chi ^{\prime }\right]  \label{2.7.4} \\
&=&\left( -1\right) ^{\deg \bar{\eta}\bar{\eta}^{\prime }}\bar{\eta}\bar{\eta%
}^{\prime }\left( \left[ \bar{\partial }\chi ,\chi ^{\prime }\right] +\left[
\chi ,\bar{\partial }\chi ^{\prime }\right] \right)  \notag \\
&=&\left[ \bar{\partial }\xi ,\xi ^{\prime }\right] +\left( -1\right) ^{\deg
\xi }\left[ \xi ,\bar{\partial }\xi ^{\prime }\right]  \notag
\end{eqnarray}

\begin{equation}
L_{\xi }\left( d\omega \right) =\left( -1\right) ^{\deg \xi }dL_{\xi }\left(
\omega \right)  \label{2.7.5}
\end{equation}

\begin{equation}
\bar{\partial }\left\langle \xi \right. \left| \omega \right\rangle
=\left\langle \bar{\partial }\xi \right. \left| \omega \right\rangle -\left(
-1\right) ^{\deg \xi }\left\langle \xi \right. \left| \bar{\partial }\omega
\right\rangle .  \label{2.7.6}
\end{equation}
And finally, from the computation 
\begin{eqnarray*}
\left\langle \left. \xi \right| L_{\xi ^{\prime }}\omega \right\rangle &=&%
\bar{\eta}\left\langle \left. \chi \right| L_{\xi ^{\prime }}\omega
\right\rangle \\
&=&\bar{\eta}\left\langle \left. \chi \right| \bar{\eta}^{\prime }L_{\chi
^{\prime }}\omega \right\rangle \\
&=&\left( -1\right) ^{\deg \bar{\eta}^{\prime }}\bar{\eta}\bar{\eta}^{\prime
}\left\langle \left. \chi \right| L_{\chi ^{\prime }}\omega \right\rangle \\
&=&\left( -1\right) ^{\deg \bar{\eta}^{\prime }}\left( \bar{\eta}\bar{\eta}
^{\prime }L_{\chi ^{\prime }}\left\langle \left. \chi \right| \omega
\right\rangle -\bar{\eta}\bar{\eta}^{\prime }\left\langle \left. \left[ \chi
^{\prime },\chi \right] \right| \omega \right\rangle \right) \\
&=&\left( -1\right) ^{\deg \xi ^{\prime }}\left\langle \left. \left[ \xi
,\xi ^{\prime }\right] \right| \omega \right\rangle +\left( -1\right)
^{\left( \deg \bar{\eta}+1\right) \deg \bar{\eta}^{\prime }}\bar{\eta}
^{\prime }L_{\chi ^{\prime }}\left\langle \left. \bar{\eta}\chi \right|
\omega \right\rangle
\end{eqnarray*}
one has 
\begin{equation}
\left( -1\right) ^{\deg \xi ^{\prime }}\left\langle \left. \xi \right|
L_{\xi ^{\prime }}\omega \right\rangle -\left( -1\right) ^{\deg \xi \deg \xi
^{\prime }}L_{\xi ^{\prime }}\left\langle \left. \xi \right| \omega
\right\rangle =\left\langle \left. \left[ \xi ,\xi ^{\prime }\right] \right|
\omega \right\rangle .  \label{2.7.7}
\end{equation}

A special case of this last formula when $\xi =\xi ^{\prime }$ and $k=1$ is 
\begin{equation}
\left[ L_{\xi },\left\langle \left. \xi \right| \ \right\rangle \right]
=\left\langle \left. \left[ \xi ,\xi \right] \right| \ \right\rangle
\label{2.7.8}
\end{equation}
and so, by type, 
\begin{equation}
\left[ L_{\xi },\left\langle \left. \xi \right| \ \right\rangle \right]
=\left\langle \left. \left[ \xi ,\xi \right] \right| \ \right\rangle =\left[
L_{\xi }^{1,0},\left\langle \left. \xi \right| \ \right\rangle \right] .
\label{2.7.9}
\end{equation}

\section{Integrability\label{A3}}

We reproduce the classical argument characterizing the systems $\left\{ \xi
\in A^{0,1}\left( T_{1,0}\left( M_{0}\right) \right) \otimes \Bbb{C}\left[
\left[ t\right] \right] \right\} $ which come from a transversely
holomorphic trivialization of a deformation $\left( \ref{2.1}\right) $.

\begin{lemma}
The ``$\left( 0,1\right) $'' tangent distribution given by the image of 
\begin{eqnarray*}
\iota :T_{0,1}\left( M_{0}\right) \rightarrow T\left( M_{0}\right) \otimes 
\Bbb{C}\left[ \left[ t\right] \right] \\
id.-\left\langle \left. \xi \right| id.\right\rangle
\end{eqnarray*}
gives, via complex conjugation, an almost complex structure on $M_{0}\times
\Delta $. This almost complex structure is integrable, that is, comes from a
(formal) deformation/trivialization of $M_{0}$ as in Definition \ref{2.2},
if and only if 
\begin{equation*}
\bar{\partial }\xi -\frac{1}{2}\left[ \xi ,\xi \right] =0.
\end{equation*}
\end{lemma}

\begin{proof}
Let $\left( v^{l}\right) $ be a system of local holomorphic coordinates on $%
M_{0}$. Locally 
\begin{equation*}
\xi =\sum\nolimits_{J}d\overline{v^{k}}\otimes h_{J,k}^{l}t^{J}\frac{
\partial }{\partial v^{l}}.
\end{equation*}
Then the image of $\iota $ is framed locally by the vector fields 
\begin{equation*}
\left( \frac{\partial }{\partial \overline{v^{k}}}-\sum\nolimits_{J,\left|
J\right| >0}h_{J,k}^{l}t^{J}\frac{\partial }{\partial v^{l}}\right)
\end{equation*}
Using a slight adaptation of the Einstein summation convention, we have 
\begin{equation*}
\begin{tabular}{l}
$\left[ \left( \frac{\partial }{\partial \overline{v^{j}}}-t^{J}h_{J,j}^{l}%
\frac{\partial }{\partial v^{l}}\right) ,\left( \frac{\partial }{\partial 
\overline{v^{k}}}-t^{J^{\prime }}h_{J^{\prime },k}^{m}\frac{\partial }{
\partial v^{m}}\right) \right] =$ \\ 
$\left[ \frac{\partial }{\partial \overline{v^{k}}},t^{J}h_{J,j}^{l}\frac{
\partial }{\partial v^{l}}\right] -\left[ \frac{\partial }{\partial 
\overline{v^{j}}},t^{J^{\prime }}h_{J^{\prime },k}^{m}\frac{\partial }{
\partial v^{m}}\right] +\left[ t^{J}h_{J,j}^{l}\frac{\partial }{\partial
v^{l}},t^{J^{\prime }}h_{J^{\prime },k}^{m}\frac{\partial }{\partial v^{m}}
\right] =$ \\ 
$t^{J}\frac{\partial h_{J,j}^{l}}{\partial \overline{v^{k}}}\frac{\partial }{
\partial v^{l}}-t^{J^{\prime }}\frac{\partial h_{J^{\prime },k}^{m}}{
\partial \overline{v^{j}}}\frac{\partial }{\partial v^{m}}+t^{J+J^{\prime
}}\left( h_{J,j}^{l}\frac{\partial h_{J^{\prime },k}^{m}}{\partial v^{l}}%
\frac{\partial }{\partial v^{m}}-h_{J^{\prime },k}^{m}\frac{\partial
h_{J,j}^{l}}{\partial v^{m}}\frac{\partial }{\partial v^{l}}\right) .$%
\end{tabular}
\end{equation*}
So integrability is checked by pairing the above vector fields with 
\begin{equation*}
dv^{r}+t^{J^{\prime \prime }}h_{J^{\prime \prime },s}^{r}d\overline{v^{s}}.
\end{equation*}
We get that integrability is equivalent to the identical vanishing of 
\begin{equation*}
t^{J}\frac{\partial h_{J,j}^{r}}{\partial \overline{v^{k}}}-t^{J^{\prime }}%
\frac{\partial h_{J^{\prime },k}^{r}}{\partial \overline{v^{j}}}
+t^{J+J^{\prime }}\left( h_{J,j}^{l}\frac{\partial h_{J^{\prime },k}^{r}}{
\partial v^{l}}-h_{J^{\prime },k}^{m}\frac{\partial h_{J,j}^{r}}{\partial
v^{m}}\right) ,
\end{equation*}
that is, of 
\begin{eqnarray*}
&&t^{J}\frac{\partial h_{J,j}^{r}}{\partial \overline{v^{k}}}d\overline{v^{k}%
}d\overline{v^{j}}-t^{J^{\prime }}\frac{\partial h_{J^{\prime },k}^{r}}{
\partial \overline{v^{j}}}d\overline{v^{k}}d\overline{v^{j}} \\
&&+t^{J+J^{\prime }}\left( h_{J,j}^{l}\frac{\partial h_{J^{\prime },k}^{r}}{
\partial v^{l}}-h_{J^{\prime },k}^{m}\frac{\partial h_{J,j}^{r}}{\partial
v^{m}}\right) d\overline{v^{k}}d\overline{v^{j}},
\end{eqnarray*}
which becomes the system of equations 
\begin{equation*}
\begin{tabular}{l}
$2d\overline{v^{k}}d\overline{v^{j}}\otimes \left( \frac{\partial h_{J,j}^{r}%
}{\partial \overline{v^{k}}}\frac{\partial }{\partial v^{r}}\right) =$ \\ 
$\sum\nolimits_{J^{\prime }+J^{\prime \prime }=J}d\overline{v^{j}}d\overline{
v^{k}}\otimes \left( h_{J^{\prime },j}^{l}\frac{\partial h_{J^{\prime \prime
},k}^{r}}{\partial v^{l}}\frac{\partial }{\partial v^{r}}-h_{J^{\prime
\prime },k}^{m}\frac{\partial h_{J^{\prime },j}^{r}}{\partial v^{m}}\frac{
\partial }{\partial v^{r}}\right) $ \\ 
$=\sum\nolimits_{J^{\prime }+J^{\prime \prime }=J}d\overline{v^{j}}d%
\overline{v^{k}}\otimes \left[ h_{J^{\prime },j}^{l}\frac{\partial }{
\partial v^{l}},h_{J^{\prime \prime },k}^{m}\frac{\partial }{\partial v^{m}}
\right] .$%
\end{tabular}
\end{equation*}
\end{proof}

\end{document}